\documentclass{article}
\usepackage{Macros}
\usepackage{hyperref}

\title{Projective Geometry for Perfectoid Spaces}
\author{Gabriel Dorfsman-Hopkins}
\date{\today}

\begin{document}
\maketitle

\begin{abstract}
To understand the structure of an algebraic variety we often embed it in various projective spaces. This develops the notion of projective geometry which has been an invaluable tool in algebraic geometry. We develop a perfectoid analog of projective geometry, and explore how equipping a perfectoid space with a map to a certain analog of projective space can be a powerful tool to understand its geometric and arithmetic structure. In particular, we show that maps from a perfectoid space $X$ to the perfectoid analog of projective space correspond to line bundles on $X$ together with some extra data, reflecting the classical theory.  Along the way we give a complete classification of vector bundles on the perfectoid unit disk, and compute the Picard group of the perfectoid analog of projective space.
\end{abstract}

\section{Introduction}
This paper is inspired by the goal of understanding vector bundles on perfectoid spaces, and how they behave under the so called \emph{tilting correspondence} of Scholze \cite{sh12}.  To do so, we develop a perfectoid analog of projective geometry.  We study the perfectoid analog of projective space defined in \cite{sh12}, which we call \emph{projectivoid space} and denote by $\bP^{n,\perf}$, and show that maps from a perfectoid space $X$ to $\bP^{n,\perf}$ correspond to line bundles on $X$ together with some extra data, giving an analog to the classical theory of maps to projective space.

To get to this point we must first understand the theory of line bundles on projectivoid space itself, and in particular, its Picard group.  In his dissertation \cite{da16}, Das worked toward computing the Picard group of the projectivoid line, $\bP^{1,\perf}$.  His proof relied on having certain local trivializations of line bundles, requiring a perfectoid analog of the Quillen-Suslin theorem.  Therefore, in order to begin developing the theory of so called projectivoid geometry, we must prove this first.

The Quillen-Suslin theorem says that finite dimensional vector bundles on affine $n$-space over a field are all trivial.  Equivalently, all finite projective modules on a polynomial ring $K[T]$ are free, where $T$ is an $n$-tuple of indeterminates.  In rigid analytic geometry, we replace polynomial rings with rings of convergent power series called \emph{Tate algebras}, denoted $K\langle T\rangle$, and it can be shown that over such rings the Quillen-Suslin theorem still holds, that is, all finite projective $K\la T\ra$-modules are free \cite{lu77}.  The analog of these rings for perfectoid spaces is the ring $K\langle T^{1/p^\infty}\rangle$ of convergent power series where the indeterminates have all their $p$th power roots.  The difficulty in extending the theorem to this \emph{perfectoid Tate algebra} is that the ring is no longer noetherian, and so the result cannot be easily reduced to the polynomial case.

In Section \ref{perfectoidtate} we study the commutative ring theoretic properties of the perfectoid Tate algebra, computing its unit group and classifying its projective modules, as well those over its integral subring and residue ring.  In particular, we establish the necessary Quillen-Suslin style result.
\begin{duplicate}[Theorem \ref{QuillenSuslinPTA} (The Quillen-Suslin Theorem for the Perfectoid Tate Algebra)]
Finite projective modules on the perfectoid Tate algebra $K\left\la T_1^{1/p^\infty},\cdots,T_n^{1/p^\infty}\right\ra$ are all free.  Equivalently, finite dimensional vector bundles on the perfectoid unit disk are all isomorphic to the trivial vector bundle.
\end{duplicate}
This completes Das' proof, and lays the groundwork to begin studying vector bundles on more general perfectoid spaces.

In Section \ref{picPn} we develop the theory of line bundles on projectivoid space, extending Das' result for $n=1$.
\begin{duplicate}[Theorem \ref{picpnperf} (The Picard Group of Projectivoid Space)]
$\Pic\bP^{n,\perf}\cong\bZ[1/p]$.
\end{duplicate}
We also compute the cohomology of all line bundles on projectivoid space, completing a computation of Bedi \cite{be18} for $n=2$.

In Section \ref{projectivoidgeometry} we compute the functor of points of projectivoid space, showing that (much like in the classical theory) it is deeply connected to the theory of line bundles on perfectoid spaces.
\begin{duplicate}[Theorem \ref{naturaliso} (The Functor of Points of Projectivoid Space)]
Let $X$ be a perfectoid space over a field $K$.  Morphisms $X\to\bP^{n,\perf}$ correspond to tuples $\left(\sL_i,s^{(i)}_j,\varphi_i\right)$, where $\sL_i\in\Pic X$, $\left\{s^{(i)}_0,\cdots,s^{(i)}_n\right\}$ are $n+1$ global sections of $\sL_i$ which generate $\sL_i$, and $\varphi_i:\sL_{i+1}^{\otimes p}\longtoo{\sim}\sL_i$ are isomorphisms under which $\left(s^{(i+1)}_j\right)^{\otimes p}\mapsto s^{(i)}_j$.
\end{duplicate}
We also provide refinements of this theorem in characteristic $p$ and see how it behaves under the tilting equivalence of Scholze.

In Section \ref{untiltingLB} we test out this new theory, using it to compare the Picard groups of a perfectoid space $X$ and its tilt $X^\tilt$.  In particular, since the tilting equivalence builds a correspondence between maps $X\to\bP^{n,\perf}_K$ and maps $X^\tilt\to\bP^{n,\perf}_{K^\tilt}$, we can chain this together with the correspondence of line bundles and maps to projectivoid space to compare line bundles on $X$ and $X^\tilt$.  The main result follows.
\begin{duplicate}[Theorem \ref{thetaInjects!}]
Suppose $X$ is a perfectoid space over $K$. Suppose that $X$ has a weakly ample line bundle (cf. Definition \ref{walb}) and that $H^0(X_{\overline K},\sO_{X_{\overline K}}) = \overline K$.  Then there is a natural injection
\[\theta:\Pic X^\tilt\into\ilim_{\sL\mapsto\sL^p}\Pic X.\]
In particular, if $\Pic X$ has no $p$ torsion, then composing with projection onto the first coordinate gives an injection
\[\theta_0:\Pic X^\flat\into\Pic X.\]
\end{duplicate}
This paper relies heavily Huber's theory of adic spaces developed in \cite{hu93} and \cite{hu94}.  We do not develop the theory here and instead refer the reader to Huber's original papers, or \cite{we12} for an excellent summary.  A summary of the theory and with an emphasis on perfectoid spaces can be found the author's doctoral dissertation (\cite{DH18} Sections 2-4), or in Kedlaya's detailed notes from the 2017 Arizona Winter School on the subject \cite{ke17}.
\subsection{Notational Conventions}
Throughout the paper we will fix a perfectoid field $K$ with a topologically nilpotent unit $\varpi$ called a \emph{pseudouniformizer}.  Its valuation ring will be denoted by $K^\circ$, whose maximal ideal is $K^{\circ\circ}$.  The residue field will be denoted by $k$.

If $R$ is a complete topological $K$-algebra we denote its subring of power bounded elements by $R^\circ$, the ideal of topologically nilpotent elements by $R^{\circ\circ}$, and the residue ring by $\tilde R$.  If $R$ is perfectoid then we denote its \emph{tilt} by $R^\tilt$.  We denote by $\sharp:R^\tilt\to R$ the Teichm\"uller map which is a map of multipicative monoids coming from composing projection onto the first coordinate with the isomorphism (of monoids) $R^\tilt\cong\ilim_{x\mapsto x^p} R$.

\subsection{Acknowledgments}
The author thanks their doctoral advisor Max Lieblich for support while pursuing this work during graduate study, and Kiran Kedlaya for his generosity of ideas while serving as host for a term at the University of California in San Diego to work on this project.  The author also thanks Peter Scholze for important comments and corrections to a previous draft.  Thanks also to Bharghav Bhatt, Daniel Bragg, Jacopo Brivio, Charles Godfrey, Tuomas Tajakka, Lucas Van Meter, Anwesh Ray, Manar Riman, Travis Scholl, Alex Voet, Peter Wear for helpful conversations.  Much of this work was completed while the author was supported under NSF grant DMS-1600813 during a semester visit to the University of California, San Diego, and under NSF grant DMS-1439786 while the author was in residence at the Institute for Computational and Experimental Research in Mathematics in Providence, RI.


\section{The Perfectoid Tate Algebra}\label{perfectoidtate}
Algebraic geometry studies the polynomial ring and it's various quotients and localizations, allowing for commutative algebraic facts to be interpreted geometrically and vice-versa. In the world of perfectoid algebras, a natural analog is the perfectoid Tate algebra (defined below), whose commutative algebra controls many of the geometric structures we study.  We therefore begin with a careful study of the algebraic object underlying most of this work.

Let $K$ be a perfectoid field with pseudouniformizer $\varpi$.  In particular $K$ is a nonarchimedean field, and so we can define the \emph{Tate algebra} $T_{n,K} = K\la X_1,\cdots,X_n\ra$ of convergent power series over $K$.  The module theory of the Tate algebra is well understood, in no small part due to  L\"utkebohmert, see \cite{lu77}.  Our main object of study in this section is the following.
\begin{Definition}[The Perfectoid Tate Algebra]\label{perftateDef}
The \emph{perfectoid Tate algebra} $T_{n,K}^\perf$ is the nonarchimedean completion
\[T_{n,K}^\perf = K\left\la X_1^{1/p^\infty},\cdots,X_n^{1/p^\infty}\right\ra = \widehat{\bigcup_{n\ge0}K\left\la X_1^{1/p^n},\cdots,X_n^{1/p^n}\right\ra}\cong\widehat{\dlim T_{n,K}}.\]
\end{Definition}
\begin{Remark}
The perfectoid Tate algebra consists of formal power series over $K$ which converge on the perfectoid unit disk.  Explicitely $X = (X_1,\cdots,X_n)$ be an $n$-tuple, we can write down the elements of this ring.
\[T_{n,K}^{\perf} = \left\{\sum_{\alpha\in\left(\bZ[1/p]_\ge0\right)^n}a_\alpha X^\alpha:\text{ for all }\lambda\in\bR_{>0}\text{ only finitely many }|a_\alpha|\ge\lambda\right\}.\]
This ring inherits the Gauss norm, $||\sum_\alpha a_\alpha X^\alpha|| = \sup\{|a_\alpha|\}$.

$\left(T_{n,K}^{\perf}\right)^\circ$ is the subring $\{||f||\le 1\}$ of power-bounded elements of $T_{n,K}^{\perf}$, and consists of power series with coefficients in $K^\circ$.  The ideal $\left(T_{n,K}^{\perf}\right)^{\circ\circ}$ of topologically nilpotent elements consist of power series with coefficients in $K^{\circ\circ}$.  The quotient is
\[\tilde T_{n,K}^{\perf}:=\left(T_{n,K}^{\perf}\right)^\circ/\left(T_{n,K}^\perf\right)^{\circ\circ} = k\left[X_1^{1/p^\infty},\cdots,X_n^{1/p^\infty}\right]\]
where $k=K^\circ/K^{\circ\circ}$ is the residue field.  Notice that every element in the quotient ring is a polynomial, because a power series in $\left(T_{n,K}^{\perf}\right)^\circ$ can only have finitely many coefficients of norm 1.
\end{Remark}
When there will be no confusion, we omit $K$ from the notation.

\begin{Remark}[A Note on Convergence]\label{convergence}
Suppose $f = \sum f_\alpha X^\alpha\in T_n^{\perf}$.  Morally speaking, saying that the sum converges should mean that evaluating $f$ at any point in the perfectoid disk should give an element of $K$.  Since are sums are not taken over $\bZ_{\ge0}^n$, but rather $(\bZ[1/p]_{\ge0})^n$, we must be more careful in defining what convergence means.  Let us begin by studying $f(1,1,\cdots,1) = \sum f_\alpha$.  This should converge, so let's begin by defining partial sums
\[s_m = \sum_{\alpha\in\left(\frac{\bZ}{p^m}\right)^n:0\le\alpha_i\le m}f_\alpha.\]
If the sequence $(s_m)$ converges, we define the infinite sum to be the limit.  Let us check that the convergence of the power series $f$ implies convergence of $\sum f_\alpha$ in this sense.  Fixing some $\varepsilon>0$, there are only finitely many $f_\alpha$ with $|f_\alpha|\ge\varepsilon$.  Therefore, there is some large $N$ such that for each such $f_\alpha$ we have $\alpha = (\alpha_1,\cdots,\alpha_n)\in\left(\frac{\bZ}{p^N}\right)^n$ $0<\alpha_i<N$.  Therefore, fixing $m\ge r>N$, the differences $s_m-s_r$ have none of the coefficients $f_\alpha$ with absolute value larger than $\varepsilon$, so that by the nonarchimedean property $|s_m-s_r|<\varepsilon$.  Thus the sum converges to an element $f(1)\in K$.

We remark now that if $|g_\alpha|\le 1$, the same argument would show that $\sum f_\alpha g_\alpha$ also converges.  This should imply that $f(x)$ converges to a point in $K$ whenever $x$ is in the perfectoid unit disk, but we defer further discussion until after we have the relevant definitions (see Remark \ref{whatAreThePoints}).
\end{Remark}
We record a useful normalization trick for further use down the line.
\begin{Lemma}[Normalization]\label{normalization}
Let $f\in T_{n}^{\perf}$ be nonzero.  There is some $\lambda\in K$ such that $||\lambda f||=1$.
\end{Lemma}
\begin{Proof}
Since only finitely many coefficients in $f$ have absolute value above $||f||-\varepsilon$, the supremum of that absolute values of the coefficients is achieved by some $f_\alpha$.  Taking $\lambda=f_\alpha^{-1}$ completes the proof.
\end{Proof}


\subsection{The Group of Units}\label{unitSection}

As a first step towards understanding the perfectoid Tate algebra, we compute its group of units.
\begin{Proposition}\label{unitz}
Let $f\in T_n^{\perf}$ with $||f||=1$.  The following are equivalent:
\begin{enumerate}[(i)]
\item{
$f$ is a unit in $\left(T_n^{\perf}\right)$.
}
\item{
$f$ is a unit in $\left(T_n^{\perf}\right)^\circ$.
}
\item{
The image of $\overline{f}$ of $f$ in $\tilde T_n^{\perf}$ is a nonzero constant $\lambda\in k^\times$.
}
\item{
$|f(0)| = 1$ and $||f-f(0)||<1$.
}
\end{enumerate}
\end{Proposition}
\begin{Proof}
(i)$\iff$(ii).  An inverse to $f$ must have absolute value 1, and therefore would also lie in $\left(T_n^{\perf}\right)^\circ$.

(ii)$\implies$(iii).  The map $\left(T_n^{\perf}\right)^\circ\to \tilde T_n^{\perf}$ must send units to units, and the group of units of $\tilde T_n^{\perf}$ is precisely the nonzero constant polynomials.  Indeed, the inverse to any element of $\tilde T_n^{\perf}$ would also have to be a polynomial (in $X^{1/p^m}$ for some $m$), implying that they both must be constants.

(iii)$\iff$(iv). This is immediate.

(iv)$\implies$(i).  If $|f(0)|=1$ then $f(0)\in K^\times\subseteq(T_n^{\perf})^\times$.  Therefore $1-\frac{f}{f(0)}\in T_n^{\perf}$ and 
\[\left|\left|1-\frac{f}{f(0)}\right|\right| = ||f(0)||\cdot\left|\left|1-\frac{f}{f(0)}\right|\right| = ||f(0)-f|| < 1.\]
Therefore $1-\frac{f}{f(0)}$ is topologically nilpotent, so that by the geometric series converges and hence $\frac{f}{f(0)}$ is a unit.  Since $f(0)$ is too, we can conclude that $f$ is a unit.
\end{Proof}
\begin{Corollary}\label{units}
$f = \sum f_\alpha X^\alpha\in\left(T_n^{\perf}\right)^\circ$ is a unit if and only if $|f_0| = 1$ and $|f_\alpha|<1$ for all $\alpha\not=0$.
\end{Corollary}
\begin{Corollary}
$f = \sum f_\alpha X^\alpha\in T_n^{\perf}$ is a unit if and only if $|f_{\alpha}|<|f_0|$ for all $\alpha\not=0$.
\end{Corollary}
\begin{Proof}
Using our normalization trick, we know $||\lambda f|| = 1$ for some $\lambda\in K^\times$.  Then $f$ is a unit if and only if $\lambda f$ is, if and only if $|\lambda f_\alpha|<1 = |\lambda f_0|$ for all $\alpha\not=0$.  Canceling shows this holds if and only if $|f_\alpha|<|f_0|$ for all $\alpha\not=0$.
\end{Proof}

\subsection{Vector Bundles on the Perfectoid Unit Disk}\label{qspta}
In classical algebraic geometry the polynomial ring (and it's various quotients) form the local building blocks of most of the objects of study.  Phrased geometrically, the prime spectrum of the polynomial ring of $n$ variables is affine space $\bA^n$ which covers or contains many of the spaces of interest.  Like schemes which are locally prime spectra of rings, perfectoid spaces are build from perfectoid algebras using the adic spectrum functor of Huber (see \cite{hu93},\cite{hu94}, and \cite{hu96}).  The perfectoid Tate algebra plays the role of the polynomial ring, and the perfectoid unit disk (defined below) plays the role of affine space.
\begin{Definition}
The \emph{perfectoid unit disk} is the adic space associated to the perfectoid Tate algebra:
$$\bD^{n,\perf} = \spa\left(T_n^{\perf},(T_n^{\perf})^\circ\right).$$
\end{Definition}
\begin{Remark}
The rigid unit disk is the adic spectrum associated to the Tate algebra. Since the perfectoid Tate algebra is the completed union of Tate algebras, using the \emph{tilde limit} formalism of \cite{shwe13}, we have
\[\bD^{n,\perf}\sim\ilim_{\varphi}\bD^{n,\ad},\]
where $\varphi$ is the $p$th power map on coordinates.  It is worth noting that the tilde limit is not the categorical inverse limit (since these are not in general unique in the category of adic spaces).  Nevertheless, it should be thought of affinoid locally as corresponding to the completed directed limit, and if such a limit exists as a perfectoid space, it is unique among all perfectoid spaces and satisfies the usual universal property (among perfectoid spaces).  See \cite{shwe13} Definition 2.4.1 and subsequent discussion.
\end{Remark}
\begin{Remark}[$K$-points of the Perfectoid Disk]\label{whatAreThePoints}
Let us compute the $K$-points of $\bD^{n,\perf}_K$, and in doing so conclude the discussion of Remark \ref{convergence}.  A $K$ point $x$ is a map of adic spaces $x:\spa(K,K^\circ)\to\bD^{n,\perf}_K$, which is equivalent to a map of Huber pairs $\varepsilon_x:\left(T_K^{n,\perf},(T_K^{n,\perf})^\circ\right)\to(K,K^\circ)$.  This is a commutative diagram of continuous $K$-algebra homomorphisms:
\begin{cd}
K\left\la X_1^{1/p^\infty},\cdots,X_n^{1/p^\infty}\right\ra\rar & K\\
K^\circ\left\la X_1^{1/p^\infty},\cdots,X_n^{1/p^\infty}\right\ra\rar\uar & K^\circ\uar.
\end{cd}
In particular, it is determined by the image of the $X_i^{1/p^k}$ in $K^\circ$.  Let $\lambda = \varepsilon_x(X_i)$.  Then $\lambda_1 = \varepsilon_x({X_i^{1/p})}$ must be a $p$th root of $\lambda$, and $\lambda_2 = \varepsilon_x(X_i^{1/p^2})$ must be a $p$th root of $\lambda_1$.  Continuing in this fashion, we see that choosing the image of the $X_i^{1/p^k}$ as $k$ varies is equivalent to fixing an element of $\ilim_{t\mapsto t^p}K^\circ = K^{\flat\circ}$.  We have therefore computed the $K$-points of the perfectoid disk:
\[\bD^{n,\perf}_K(K) = (K^{\flat\circ})^n.\]
The evaluation function $\varepsilon_x(f)=f(x)$ then amounts to plugging the coordinates of $x$ as an element of $(K^{\flat\circ})^n$ into $f$, which we saw in Remark \ref{convergence} converges to a point in $K$.
\end{Remark}

There is a well known correspondence between finite projective modules over a ring, and finite dimensional (algebraic) vector bundles over the associated affine scheme, and more generally, between vector bundles over a locally ringed space and locally free sheaves on that space (see, for example, \cite{ha77} exercise 2.5.18).  In \cite{se55}, Serre conjectured that all finite projective modules over the polynomial ring $A = k[x_1,\cdots,x_n]$ are free.  This can be interpreted geometrically as saying there are no nontrivial algebraic vector bundles over affine space $\bA^n = \spec A$.  In 1976, Quillen \cite{qu76} and independently Suslin proved Serre's conjecture, which is now known as the Quillen-Suslin theorem.  L\"utkebohmert in \cite{lu77} was shortly able to extend the result to the Tate algebra $K\la X_1,\cdots,X_n\ra$ of convergent power series over a complete nonarchimedean field.

In what follows, we prove a perfectoid analog of the Quillen-Suslin theorem. Specifically, we prove that all finite projective modules on the perfectoid Tate algebra are trivial.  This will imply that the perfectoid unit disk has no nontrivial finite vector bundles.  Along the way we will show that both the subring of integral elements $\left(T_n^\perf\right)^\circ$, and the residue ring $\tilde T_n^\perf$ also have no nontrivial finite projective modules.  Although these results are not necessary to establish the result for the perfectoid Tate algebra, they will be important in asserting the effectiveness of the \v Cech cohomology groups of certain sheaves in Section \ref{picPn}.

\subsubsection{Finite Projective Modules on the Residue Ring}\label{residuesection}

Let us begin by proving that finite projective modules are free over the residue ring
\[\tilde T_n^{\perf} = k\left[X_1^{1/p^\infty},\cdots,X_n^{1/p^\infty}\right] = \bigcup_m k\left[X_1^{1/p^m},\cdots,X_n^{1/p^m}\right].\]
To see this, we first briefly review a (non unique) correspondence between projective modules and idempotent matrices.

Let $R$ be a commutative ring, and fix a finite projective $R$ module $P$. Consider a presentation $\pi:R^n\to P$, as well as a section of this projection $\sigma$.  The composition $\sigma\circ\pi$ produces an idempotent matrix $U\in M_n(R)$.
\begin{cd}
R^n\ar["\pi"]{dr}\ar["U"]{rr}\ar["U",bend left=30]{rrrr} & {} & R^n\ar["\pi"]{dr}\ar["U"]{rr} & {} & R^n\\
{} & P\ar["\sigma"]{ur}\ar[equals]{rr} & {} & P\ar["\sigma"]{ur} &{}
\end{cd}
Conversely, the image of an idempotent matrix $U$ is always projective, with the section just given by the natural inclusion $\im U\subseteq R^n$.  In this way, we get a (non unique) correspondence between finite projective modules and idempotent matrices over $R$.
\begin{Lemma}\label{colimits}
Suppose $R = \dlim_i R_i$ is a filtered colimit of commutative rings.  Then every finite projective $R$ module is the base extension of a finite projective $R_i$-module.
\end{Lemma}
\begin{Proof}
To a finite projective $R$ module $M$ we may associate a projector matrix $U$. Each entry in the matrix is defined over some $R_i$, and as the limit is filtered, $U$ is defined over $R_i$ for some (perhaps larger) $i$.  Its image as a map from $R_i^n$ to itself is therefore a projective $R_i$-module whose base extension is $M$.
\end{Proof}

\begin{Corollary}\label{reductionqs}
Let $R = k\left[X_1^{1/p^\infty},\cdots,X_r^{1/p^\infty},X_{r+1}^{\pm 1/p^\infty},\cdots,X_n^{\pm 1/p^\infty}\right].$  Every $R$-module is free.  In particular (letting $r=n$), every $\tilde T_n^{\perf}$-module is free.
\end{Corollary}
\begin{Proof}
As $R$ is the filtered colimit of (Laurent) polynomial rings, Lemma \ref{colimits} implies that a finite projective $R$-module $M$ is the base extension of some $N$ over a (Laurent) polynomial ring.  By the Quillen-Suslin theorem, $N$ is free, so $M$ is too.
\end{Proof}

\subsubsection{Finite Projective Modules on the Subring of Integral Elements}

We extend Corolary \ref{reductionqs} to the subring of power-bounded elements of the perfectoid Tate algebra, $\left(T_n^{\perf}\right)^\circ$, using Nakayama's lemma.  We first fix some notation.
\begin{Notation}
For a commutative ring $R$ and an ideal $I$ contained in the Jacobson radical of $R$ we let $R_0 = R/I$.  For an $R$ module $M$ we will denote by $M_0$ the $R_0$ module $M/IM$, and for a homomorphism $\phi$ of $R$-modules we denote by $\phi_0$ its reduction mod $I$.  If $m\in M$, then we denote by $\overline m$ its image in $M_0$.
\end{Notation}
\begin{Lemma}\label{reduceit}
Let $R$ be a commutative ring, and $I$ an ideal contained in the Jacobson radical of $R$.  If $M$ and $N$ are two projective $R$-modules such that there exists an isomorphism $\phi: M_0\too{\sim} N_0$, then $\phi$ lifts to an isomorphism $\psi:M\too{\sim} N$.
\end{Lemma}
\begin{Proof}
We have the following commutative diagram.
\begin{cd}
M\ar[d, two heads, "\pi"]\ar[r, dashed, "\psi"] & N\ar[d,two heads,"\rho"]\\
M_0\ar[r,"\phi"] & N_0
\end{cd}
Indeed, a lift $\psi$ exists because $M$ is projective.  Notice that $\psi_0 = \phi$.  Indeed,
\begin{eqnarray*}
\psi_0(\overline m) &=& \overline{\psi(m)}\\
&=& \rho\psi(m)\\
&=&\phi\pi(m)\\
&=&\phi(\overline m).
\end{eqnarray*}
Since $\phi$ surjects, so does $\psi$ by Nakayama's lemma.  Since $N$ is projective, $\psi$ has a section $\sigma$ which is necessarily injective.  We claim that $\sigma_0 = \phi^{-1}$.  Indeed, we can check this after applying $\phi$.
\begin{eqnarray*}
\phi\sigma_0(\overline{n}) &=& \phi\pi\sigma(n)\\
&=&\rho\psi\sigma(n)\\
&=&\rho(n)\\
&=&\overline n.
\end{eqnarray*}
Therefore $\sigma_0$ surjects, so that $\sigma$ surjects by Nakayama's lemma.  Thus $\sigma$ is an isomorphism.
\end{Proof}
\begin{Corollary}\label{projtofree}
With the same setup as Lemma \ref{reduceit}, we let $P$ be a projective $R$-module.  If $P_0$ is a free $R_0$ module, then $P$ is free.
\end{Corollary}
\begin{Proof}
Suppose $P_0\cong R_0^m$.  We also have $(R^m)_0\cong R_0^m$ so that by Lemma \ref{reduceit}, $P\cong R^m$.
\end{Proof}
We now apply this to the case at hand.
\begin{Lemma}\label{jacobsyeh}
Let $R$ be one of the following:
\begin{itemize}
\item{
$\left(T_n^{\perf}\right)^\circ$ or one of its associated Laurent series rings given by inverting a coordinate function.
}
\item{Either of the previous two rings modulo a positive power of $\varpi$.}
\end{itemize}
Explicitly (and up to a rearranging of coordinates),
\[R = K^\circ\left\la X_1^{1/p^\infty},\cdots,X_r^{1/p^\infty},X_{r+1}^{\pm 1/p^\infty},\cdots,X_n^{\pm 1/p^\infty}\right\ra,\]
or
\[R = K^\circ/\varpi^d \left[ X_1^{1/p^\infty},\cdots,X_r^{1/p^\infty},X_{r+1}^{\pm 1/p^\infty},\cdots,X_n^{\pm 1/p^\infty}\right].\]
The ideal of topologically nilpotent (resp. nilpotent) elements lies in the Jacobson radical of $R$.
\end{Lemma}
\begin{Proof}
If $f$ is (topologically) nilpotent then so is $fg$ for all $g$.  Thus the geometric series for $\frac{1}{1-fg}$ converges to an inverse of $1-fg$, so that it is a unit.  Since $g$ was arbitrary, this shows that $f$ is in the Jacobson radical.
\end{Proof}
The desired result follows.
\begin{Corollary}\label{qsbounded}
Let $R$ as in Lemma \ref{jacobsyeh}.  Every finite projective $R$-module is free.
\end{Corollary}
\begin{Proof}
Notice that $R_0$, the reduction of $R$ modulo the ideal of (topologically) nilpotent elements is the ring of Corollary \ref{reductionqs}. Let $P$ be a finite projective $R$ module.  Then $P_0$ is a finite projective $R_0$ module and therefore is free.  Therefore by Corollary \ref{projtofree} it suffices to show that the kernel of the reduction map is contained in the Jacobson radical, but this is Lemma \ref{jacobsyeh}.
\end{Proof}


\subsubsection{The Quillen-Suslin Theorem for the Perfectoid Tate Algebra}
We now prove the main result of this section.
\begin{Theorem}\label{QuillenSuslinPTA}
Finite projective modules on the perfectoid Tate algebra are free.  Equivalently, all finite vector bundles on the perfectoid unit disk are free.
\end{Theorem}
In fact, this will follow from something slightly more general.
\begin{Theorem}\label{laurents}
Let $R$ be a perfectoid Laurent series algebra:
\[R = K\left\la X_1^{1/p^\infty},\cdots,X_r^{1/p^\infty},X_{r+1}^{\pm 1/p^\infty},\cdots,X_n^{\pm 1/p^\infty}\right\ra.\]
Then every finite projective $R$ module is free.
\end{Theorem}
A result of Gabber and Romero will do most of the heavy lifting, and state the result here:
\begin{Proposition}[\cite{gr03} Corollary 5.4.42]\label{gabber}
Let $R$ be a commutative ring, $t\in R$ a nonzero divisor, and $I\subset R$ an ideal.  Let $\hat R$ be the $tI$-adic completion of $R$, and suppose $(R,tI)$ form a Henselian pair.  Then the base extension functor $R[t^{-1}]-\textbf{Mod}\to \hat R[t^{-1}]-\textbf{Mod}$ induces a bijection between isomorphism classes of finite projective $R[t^{-1}]$-modules and finite projective $\hat R[t^{-1}]$-modules.
\end{Proposition}
Let us describe the objects we will feed into this result.  Playing the role of $R$ is:
\begin{equation}\label{defR}
R = \bigcup_i K^\circ\left\la X_1^{1/p^i},\cdots,X_r^{1/p^i},X_{r+1}^{\pm 1/p^i},\cdots,X_n^{\pm 1/p^i}\right\ra =: \bigcup_i R_i.\end{equation}
Playing the role of $t$ we have $\varpi$, and $I$ will be the unit ideal.  Therefore
\[R[1/\varpi]=\bigcup_i K\left\la X_1^{1/p^i},\cdots,X_r^{1/p^i},X_{r+1}^{\pm 1/p^i},\cdots,X_n^{\pm 1/p^i}\right\ra.\]
Then $\hat R$ is the ring of integral elements (cf. Lemma \ref{jacobsyeh}), and $\hat R[1/\varpi]$ is the perfectoid Tate algebra (or Laurent series algebra), for which we are trying to prove the Quillen-Suslin result.  Theorem \ref{laurents} will follow from the following two lemmas.
\begin{Lemma}\label{reduceToTate}
Finite projective modules on $R[1/\varpi]$ are all free.
\end{Lemma}
\begin{Proof}
Applying Lemma \ref{colimits}, a finite projective $R[1/\varpi]$-module $M$ is the base extention of a finite projective module $N$ over a rigid analytic Tate algebra (or Laurent series algebra).  But $N$ is free due to the rigid analytic Quillen-Suslin theorem (\cite{lu77} Satz 1), so we are done.
\end{Proof}
\begin{Lemma}\label{hensel}
The pair $(R,(\varpi))$ is Henselian.
\end{Lemma}
\begin{Proof}
Lemma \ref{jacobsyeh} shows that $\varpi$ is in the Jacobson radical of $R$.  Suppose $f\in R$ is monic, and that $\overline f = g_0h_0\in R/\varpi$ with $g_0,h_0$ monic.  In fact, $f\in R_i$ for some $i$ (using the notation of Equation \ref{defR} above).  Perhaps increasing $i$, we can take the factorization of $f$ to take place in $R_i/\varpi$.  As $R_i$ is $\varpi$-complete, $(R_i,\varpi)$ form a Henselian pair (see, for example, \cite{stacks-project} Tag 0ALJ).  Therefore the factorization lifts to $f=gh$ in $R_i\subset R$, with $g$ and $h$ monic.
\end{Proof}
The main result now follows easily.
\begin{nProof}[Proof of Theorem \ref{laurents}]
Lemma \ref{hensel} allows us to apply Proposition \ref{gabber}, and conclude that each finite projective module over $\hat R[1/\varpi]$ is the base extension of one on $R[1/\varpi]$, which by Lemma \ref{reduceToTate} must be free.
\end{nProof}

\section{Line Bundles and Cohomology on Projectivoid Space}\label{picPn}

In classical algebraic geometry, the notion of \emph{projective geometry} is a very powerful tool to study properties of varieties and schemes.  Indeed, one can learn a lot about a scheme by understanding its maps to various projective spaces, and this theory is intimately connected to the theory of line bundles on that space.  In this and the following section we develop an analogous theory for perfectoid spaces.  Let us begin by defining a perfectoid analog of projective space.

\begin{Definition}[The Projectivoid Line]\label{projline}
Analogously to the construction of the Riemann sphere, we can build the perfectoid analog of the projective line by gluing two copies of the perfectoid unit disk along the perfectoid unit circle.

Explicitly, the inclusion $K\left\la T^{1/p^\infty}\right\ra\to K\left\la T^{\pm1/p^\infty}\right\ra$ corresponds to the open immersion $\bS^{1,\perf}\into\bD^{1,\perf}$.  The map $K\left\la T^{-1/p^\infty}\right\ra\to K\left\la T^{\pm1/p^\infty}\right\ra$ also corresponds to an open immersion of the circle into the disk, where now the disk has coordinate $T^{-1}$.  Identifying the circles on each of these disks and gluing produces \emph{the projectivoid line}, denoted $\bP^{1,\perf}_K$.
\end{Definition}
\begin{Definition}[Projectivoid Space]\label{projectivoidspace}
As with the projectivoid line, we can define projectivoid $n$-space by gluing together $n+1$ perfectoid unit $n$-polydisks along their associated perfectoid sphere exteriors as in Example \ref{projline}.  This equips projectivoid space with a cover by perfectoid unit disks which we will henceforth refer to as the \emph{standard cover}.

In \cite{sh14} Section 7, Scholze showed that we could also define projectivoid space in the following way.  Let $\bP_K^n$ be projective space over $K$, which can be viewed as an adic space as in by first being viewed as a rigid space using the rigid analytification functor, and then as an adic space as in \cite{hu94}.  Let $\varphi:\bP_K^n\to\bP_K^n$ be the morphism given in projective coordinates by $(T_0:\cdots:T_n)\mapsto(T_0^p:\cdots:T_n^p)$.  Then
\[\bP_K^{n,\perf}\sim\ilim_{\varphi}\bP_K^n.\]
As with the perfectoid disk ``$\sim\ilim$'' is the `tilde limit' of \cite{shwe13}.
\end{Definition}
\begin{Remark}
Notice that all the finite intersections of the standard cover will be the adic spectra associated to what we called perfectoid Laurent series algebras in the previous section.  In particular, Theorem \ref{laurents} shows that all their finite vector bundles are free.
\end{Remark}

Scholze showed in \cite{sh12} that the construction of projectivoid space is compatible with the tilting functor.  In this section we begin our exploration of so called \emph{projectivoid geometry} by developing the theory of line bundles on projectivoid space.  In particular, we compute the Picard group of $\bP^{n,\perf}$, as well as the sheaf cohomology of all line bundles.  We continue developing the theory in the following section, where we will show how an arbitrary perfectoid space's maps to projectivoid space is intimately connected to its theory of line bundles, reflecting the situation in classical algebraic geometry, but with an extra arithmetic twist. 
\subsection{The Picard Group of Projectivoid Space}\label{computingpic}
The main result of this section follows:
\begin{Theorem}\label{picpnperf}
$\Pic\bP^{n,\perf}\cong\bZ[1/p]$.
\end{Theorem}
We begin by outlining our proof.  We use that for any ringed space $X$, there is a natural isomorphism $\Pic(X)\cong\HH^1(X,\sO_X^*)$ (see for example \cite{ha77} Exercise III.4.5).  Let $\sO_X^+\subset\sO_X$ be the subsheaf of integral functions, defined by the rule
\[U\mapsto\{f\in\sO_X(U):|f(x)|\le1\text{ for all }x\in U\}.\]
It is a sheaf of rings, and we let $\sO_X^{+*}$ be it's unit group.  We begin in Section \ref{cechreductions} by establishing an isomorphism
\[\HH^1(X,\sO_X^*)\longtoo{\sim}\HH^1(X,\sO_X^{+*}).\]
Let $\sO_X^{++}\subset\sO_X^+$ be the sheaf of topologically nilpotent functions, defined by the rule
\[\sO_X^{++}:U\mapsto\{f\in\sO_X(U):|f(x)|<1\text{ for all }x\in U\}.\]
It is a sheaf of ideals in $\sO_X^+$ and the quotient we denote by $\tilde \sO_X$.  The map $\sO_X^{+*}\to\tilde\sO_X^*$ induces a map on cohomology, $\HH^1(X,\sO_X^{+*})\longto\HH^1(X,\tilde\sO_X^*)$.  The latter is easier to understand, and we finish Section \ref{cechreductions} by directly computing that it is isomorphic to $\bZ[1/p]$.  Composing all this gives us a map
\[\varphi:\Pic X\longto\bZ[1/p].\]
In Section \ref{deformation} we show that $\varphi$ is an isomorphism by deforming invertible modules along the sheaves of algebras $\sA_d = \sO_X^+/\varpi^d$, which live between $\sO_X^+$ and $\tilde\sO_X$.

We make frequent use of the fact projectivoid space comes equipped with a standard cover by perfectoid disks whose geometry we understand well due to the results of Section \ref{perfectoidtate}.  In particular, any line bundle on $\bP^{n,\perf}$ becomes trivial on the standard cover and its various finite intersections due to Theorem \ref{laurents}.  We therefore use \v{C}ech cohomology with respect to this cover to study line bundles on projectivoid space.  The passage between \v Cech and sheaf cohomology is safe due to the following lemma.
\begin{Lemma}\label{cechOK}
There are natural isomorphisms
\begin{itemize}
\item $\HH^1(X,\sO_X^*) \cong \cHH^1(\fU,\sO_X^*)$.
\item $\HH^1(X,\sO_X^{+*}) \cong \cHH^1(\fU,\sO_X^{+*})$.
\item $\HH^1(X,\sA_d^*) \cong \cHH^1(\fU,\sA_d^*)$ for all $d$.
\item $\HH^1(X,\tilde\sO_X^{*}) \cong \cHH^1(\fU,\tilde\sO_X^{*})$.
\end{itemize}
\end{Lemma}
\begin{Proof}
Let $\sR$ be one of the following sheaves of rings: $\sO_X,\sO_X^+,\sA_d$, or $\tilde\sO_X$, and let $\sR^*$ be the associated sheaf of units.  The \v Cech to Derived spectral sequence (\cite{SGA4} Expos\'e V Th\'eor\`eme  3.2) describes the following spectral sequence:
\[E_2^{p,q}:\cHH^p(\fU,\sH^q(\sR^*)) \implies \HH^{p+q}(X,\sR^*).\]
Let $U$ be a finite intersection of elements in the standard cover.  Then either Theorem \ref{laurents}, Corollary \ref{qsbounded}, or Corollary \ref{reductionqs} implies that $\sR(U)$ has no nontrivial invertible modules.  In particular $\sH^1(\sR^*)(U) = H^1(U,\sR^*)=0$, as such a cohomology class would construct an invertible $\sR(U)$ module.  In particular, the sequence of low degree terms for the spectral sequence degenerates to:
\[\cHH^1(\fU,\sR^*)\cong\HH^1(X,\sR^*).\]
\end{Proof}

\subsubsection{Reduction to the Cohomology of Unit Groups of the Integral and Residue Sheaves}\label{cechreductions}
Lemma \ref{cechOK} implies that $\Pic(X)\cong\cHH^1(\fU,\sO_X^*)$.  Our first step is to compare this to the cohomology on the sheaf of integral units.
\begin{Lemma}\label{integralunitcoho}
$\Pic(X)\cong\cHH^1(\fU,\sO_X^{+*})$.
\end{Lemma}
\begin{Proof}
We have a short exact sequence of chain complexes:
\begin{cd}
0 \rar & \prod_i \sO_X^{+*}(U_i)\rar\dar & \prod_{i}\sO_X^*(U_i)\rar\dar & \prod_{i}|K^*|\rar\dar & 0\\
0 \rar & \prod_{i,j} \sO_X^{+*}(U_i\cap U_j)\rar\dar & \prod_{i,j}\sO_X^*(U_i\cap U_j)\rar\dar & \prod_{i,j}|K^*|\rar\dar & 0\\
0 \rar & \prod_{i,j,k} \sO_X^{+*}(U_i\cap U_j\cap U_k)\rar\dar & \prod_{i,j,k}\sO_X^*(U_i\cap U_j\cap U_k)\rar\dar & \prod_{i,j,k}|K^*|\rar\dar & 0\\
{}& \vdots & \vdots & \vdots & {}
\end{cd}
The left and middle complexes are the \v{C}ech complexes for $\sO_X^{+*}$ and $\sO_X^*$ respectively, and the map on the right is $|\cdot|$ which is plainly surjective.  Also the right hand complex has kernel $|K^*|$ and is otherwise exact, so that the long exact sequence on cohomology gives
\[\cHH^i(\fU,\sO_X^{+*})\cong\cHH^i(\fU,\sO_X^*),\]
for all $i>0$.  Letting $i=1$ completes the proof.
\end{Proof}
With this in hand, we can study what happens to cohomology classes reducing modulo the sheaf of topologically nilpotent functions.  Indeed, the surjection $\sO_X^+\to\tilde\sO_X$ induces a map on unit groups $\sO_X^{+*}\to\tilde\sO_X^*$.  Taking cohomology gives
\[\HH^1(X,\sO_X^{+*})\to\HH^1(X,\tilde\sO_X^*).\]
The source we know is $\Pic X$ by Lemmas \ref{cechOK} and \ref{integralunitcoho}.  We next compute the target.
\begin{Lemma}
$\HH^1(X,\tilde\sO_X^*)\cong\bZ[1/p]$.
\end{Lemma}
\begin{Proof}
By Lemma \ref{cechOK} it suffices to compute $\cHH^1(\fU,\tilde\sO_X^*)$.  We will make use of projective coordinates $[T_0:\cdots:T_n]$ for $\bP^{n,\perf}$, and as above we denote by $k$ the residue field of $K$.  We denote the differentials of the \v{C}ech complex $\cech^*(\fU,\tilde\sO_X)$ by $d^i$. Note that for all $U_i$ we have
\[\tilde\sO_X(U_i) = k\left[\left(\frac{T_0}{T_i}\right)^{1/p^\infty},\cdots,\left(\frac{T_n}{T_i}\right)^{1/p^\infty}\right],\]
so that $\tilde\sO_X^*(U_i)\cong k^*$, since the only invertible polynomials are the constant functions.  Therefore we have $C^0(\fU,\tilde\sO_X^*(U_i))\cong (k^*)^{n+1}$, and viewing the kernel as the intersection we have $\ker d^0\cong k^*$, so that $\im d^0\cong (k^*)^n$.

The rings $\tilde\sO_X(U_i\cap U_j)$ consist of Laurent polynomials, and the only invertible Laurent polynomials are monomials in the invertible variable.  That is,
\[\tilde\sO_X^*(U_i\cap U_j) = \left\{\lambda\left(\frac{T_i}{T_j}\right)^\alpha : \lambda\in k^*,\alpha\in\bZ[1/p]\right\}\cong k^*\oplus\bZ[1/p].\]
Let $(f_{ij})\in C^1(\fU,\tilde\sO_X^*)$, and suppose $(f_{ij})\in\ker d^1$.  This means that for all $i<j<k$ we have $f_{ij}f_{jk} = f_{ik}$.  That is,
\[\lambda_{ij}\left(\frac{T_i}{T_j}\right)^{\alpha_{ij}}\cdot\lambda_{jk}\left(\frac{T_j}{T_k}\right)^{\alpha_{jk}} = \lambda_{ik}\left(\frac{T_i}{T_k}\right)^{\alpha_{ik}}.\]
In particular, $\alpha_{ij} = \alpha_{jk} = \alpha_{ik}$, and so the degree of every factor in an element of the kernel must match.  The fact that $\lambda_{ij}\lambda_{jk}=\lambda_{ik}$ leaves $n$ degrees of freedom for the coefficient, so that $\ker d^1 = (k^*)^n\oplus\bZ[1/p]$, and since $\im d^0 = (k^*)^n$, we conclude that $\cHH^1(\fU,\tilde\sO_X)\cong\bZ[1/p]$. 
\end{Proof}
Making the necessary identifications gives a map $\varphi:\Pic X\to\bZ[1/p]$.  It remains to prove it is an isomorphism.  For now we show it is surjective.
\begin{Lemma}\label{sectionexists}
The map $\varphi:\Pic X\to\bZ[1/p]$ has a natural section.
\end{Lemma}
\begin{Proof}
The $\bZ[1/p]$ on the right can be interpreted as follows.  For each $\alpha\in\bZ[1/p]$ we can build an invertible $\tilde\sO_X$-module starting with $\tilde\sO_{U_i}$ on each open set in the standard cover, and gluing along transition maps $\left(\frac{x_i}{x_j}\right)^\alpha$ on $U_i\cap U_j$.  Doing the same construction, but starting with $\sO_{U_i}^+$ gives us a section of $\varphi$ (call it $\sigma$).  In particular, $\varphi$ is surjective and we have an embedding $\bZ[1/p]\into\Pic(X)$.
\end{Proof}
\begin{Remark}
Therefore we can consider twisting sheaves $\sO(d)\in\Pic(X)$ for every $d\in\bZ[1/p]$.
\end{Remark}
We also have the following immediate corollary.
\begin{Corollary}
Let $A = K\left\la T_0^{1/p^\infty},\cdots,T_n^{1/p^\infty}\right\ra$.  Then,
\[\Gamma(U_{i_1}\cap\cdots\cap U_{i_r},\sO(d))\cong\left(\widehat{A_{T_{i_1}\cdots T_{i_r}}}\right)_d\]
the degree $d$ part of the completion of the localization of $A$.
\end{Corollary}
\subsubsection{Deforming Line Bundles from the Residue}\label{deformation}
Most of the following computation will be made using \v Cech cohomology, with $\fU$ the standard cover of $X=\bP^{n,\perf}$ as above.  For the crucial deformation steps we will need the following lemma, whose proof we defer to the following section (see Remark \ref{integralcohomology}).
\begin{Lemma}\label{cechacyclic}
$\cHH^i(\fU,\sO_X^+)=0$ for all $i>0$.
\end{Lemma}
We will also need the following lemma of commutative algebra.
\begin{Lemma}\label{Jacobson}
Let $R\to S$ be a surjection of rings whose kernel $I$ is contained in the Jacobson radical of $R$.  Then the induced map on unit groups, $R^*\to S^*$, remains surjective.
\end{Lemma}
\begin{Proof}
Fix $s\in S^*$ and $r\in R$ mapping to $s$.  If $r\in\fm$ for any maximal ideal of $r$, then its image would be contained in $\fm/I\cdot\fm$, a proper ideal of $S$.  Since its image is a unit this cannot be the case.  Since $r$ is not contained in any maximal ideal it must be a unit.
\end{Proof}
Let us begin by enumerating a few useful exact sequences.
\begin{equation}\label{*1}
0\longto\sO_X^{++}\longto\sO_X^+\longto\tilde\sO_X\longto 0.
\end{equation}
Because $\sO_X^{++}$ consists of toplogically nilpotent functions, it is contained in the Jacobson radical of $\sO_X^+$, so that by Lemma \ref{Jacobson} the right hand map of the sequence remains surjective on unit groups.
\begin{equation}\label{*2}
1\longto 1+\sO_X^{++}\longto \sO_X^{+*}\longto\tilde\sO_X^*\longto1.
\end{equation}
If $1+\sO_X^{++}$ were acyclic, we could reduce finding line bundles on $X$ to finding invertible $\tilde\sO_X$ modules and be done, but this acyclicity has so far been rather elusive (and seems unlikely in general).  That being said, there is a filtration of $\sO_X^{++}$ by sheaves of (principal) ideals $\varpi^d$ for $d>0$.
\begin{equation}\label{*3}
0\longto\varpi^d\longto\sO_X^+\longto\sA_d\longto0.
\end{equation}
Notice that $\sA_d$ is a sheaf of (nonreduced) $\sO_X^+$-algebras, and that for every $d'>d$, we have surjections $\sA_{d'}\onto\sA_d$ with kernel $\varpi^d/\varpi^{d'}$.  As before, $\varpi^d$ is contained in the Jacobson radical, so that we also have
\begin{equation}\label{*4}
1\longto 1+\varpi^d\longto\sO_X^{+*}\longto\sA_d^*\longto1.
\end{equation}
\begin{Lemma}\label{limitsAndColimits}
The sheaves $\sA_d^*$ interpolate continuously between $\tilde\sO_X^*$ and $\sO_X^{+*}$.  More precisely:
\begin{itemize}
\item $\dlim\sA_d^*\cong\tilde\sO_X^*$
\item $\ilim\sA_d^*\cong\sO_X^{+*}.$
\end{itemize}
\end{Lemma}
\begin{Proof}
For the first statement notice that $\dlim(1+\varpi^d) = \cup_d(1+\varpi^d) = 1 +\sO_X^{++}$.  This follows if $\sO_X^{++} = \cup_d(\varpi^d)$, so fix $f\in\sO_X^{++}$ so that have $|f|=\lambda<1$.  There is some $d$ such that $1>|\varpi^d|>\lambda$, so that $|f/\varpi^d|<1$.  Thus $f'=f/\varpi^d\in\sO_X^++$ and $f\in\varpi^d$.  Since colimits of abelian sheaves are exact, applying $\dlim_d$ to Sequence \ref{*3} proves part 1 of the lemma.

For part 2, notice that $\ilim(\sA_d)\cong\sO_X^+$ since $\sO_X^+$ is $\varpi$-adically complete.  Since the unit group functor commutes with inverse limits (indeed, it is left adjoint to the group ring functor), are done.
\end{Proof}
We next assert that Sequences \ref{*1} through \ref{*4} induce long exact sequences on \v Cech cohomology.
\begin{Lemma}\label{longexactcech}
Each of the sequences \ref{*1} through \ref{*4} remain exact when evaluated on any finite intersection of elements in the standard cover $\fU$.  In particular, they each induce long exact sequences on \v Cech cohomology.
\end{Lemma}
\begin{Proof}
We need only check exactness on the right.  For sequences \ref{*1} and \ref{*3} these evaluate to:
\[K^\circ\left\la X_1^{1/p^\infty},\cdots,X_r^{1/p^\infty},X_{r+1}^{\pm 1/p^\infty},\cdots,X_n^{\pm 1/p^\infty}\right\ra\longto k\left[ X_1^{1/p^\infty},\cdots,X_r^{1/p^\infty},X_{r+1}^{\pm 1/p^\infty},\cdots,X_n^{\pm 1/p^\infty}\right]\]
and
\[K^\circ\left\la X_1^{1/p^\infty},\cdots,X_r^{1/p^\infty},X_{r+1}^{\pm 1/p^\infty},\cdots,X_n^{\pm 1/p^\infty}\right\ra\longto K^\circ/\varpi^d\left[X_1^{1/p^\infty},\cdots,X_r^{1/p^\infty},X_{r+1}^{\pm 1/p^\infty},\cdots,X_n^{\pm 1/p^\infty}\right]\]
which are plainly surjective.  Sequences \ref{*2} and \ref{*4} come from applying the unit group functor to surjections above.  Lemma \ref{jacobsyeh} implies that the kernel of each is contained in the Jacobson radical, so that the maps remain surjective on unit groups by Lemma \ref{Jacobson}.
\end{Proof}

\begin{Lemma}\label{l1}
For all $i>0$ and $d>0$, $\cHH^i(\fU,\sA_d)=0$.
\end{Lemma}
\begin{Proof}
This follows from the long exact sequence on \v Cech cohomology associated to Sequence \ref{*3} and Lemma \ref{cechacyclic}, noticing that $\varpi^d\cong\sO_X^+$ since it is a principal ideal.
\end{Proof}

\begin{Lemma}\label{l2}
For all $d,i>0$, the natural map
\[\cHH^i(\fU,\sA_{2d}^*)\longto\cHH^i(\fU,\sA_d^*)\]
is an isomorphism.
\end{Lemma}
\begin{Proof}
Consider:
\begin{cd}
0\rar & 1+\varpi^{2d}\rar\ar[hookrightarrow]{d} & \sO_X^{+*} \ar[equal]{d} \rar &\sA_{2d}^*\rar\ar[twoheadrightarrow]{d} & 0\\
0\rar & 1 + \varpi^d\rar & \sO_X^{+*}\rar & \sA_d^*\rar & 0
\end{cd}
By the snake lemma, we have
\begin{equation}\label{*5}
1\longto1+\varpi^d/\varpi^{2d}\longto\sA_{2d}^*\longto\sA_d^*\longto1.
\end{equation}
Notice that $1+\varpi^d/\varpi^{2d}$ has a natural $\sA_d$-module structure making it isomorphic to $\sA_d$, given locally by the map $1\mapsto\varpi^d$.  Indeed, the map is well defined because $\varpi^d/\varpi^{2d}$ is a square zero ideal, and the kernel is precisely $\varpi^d$ (which is 0 in $\sA_d$), while surjectivity is clear.  In particular, by Lemma \ref{l1}, $1+\varpi^d/\varpi^{2d}$ has no higher \v Cech cohomology, so the conclusion follows if there is a long exact sequence on \v Cech cohomology associated Sequence \ref{*5}.

Let $U$ be a finite intersection of elements of the standard cover.  It remains to show that Sequence \ref{*5} remains exact when evaluated on $U$.  By Lemma \ref{longexactcech},
\begin{cd}
0\rar & \Gamma(U,1+\varpi^{2d})\rar\ar[hookrightarrow]{d} & \Gamma(U,\sO_X^{+*}) \ar[equal]{d} \rar &\Gamma(U,\sA_{2d}^*)\rar\ar{d} & 0\\
0\rar & \Gamma(U,1 + \varpi^d)\rar & \Gamma(U,\sO_X^{+*})\rar & \Gamma(U,\sA_d^*)\rar & 0
\end{cd}
has exact rows.  The ring map $\sA_{2d}(U)\to\sA_d(U)$ is surjective with nilpotent kernel, so that again applying Lemma \ref{Jacobson} we see the vertical map on the right is surjective.  Thus the snake lemma exhibits Sequence \ref{*5} evaluated at $U$ as an exact sequence, completing the proof.
\end{Proof}

\begin{Lemma}\label{l3}
For all $d'>d>0$ and $i>0$, the natural map:
\[\cHH^i(\fU,\sA_{d'}^*)\longto\cHH^i(\fU,\sA_d^*),\]
is an isomorphism.
\end{Lemma}
\begin{Proof}
By Lemma \ref{l2}, replacing $d$ with $2^ld$ will not change \v Cech cohomology, so we may assume $d<d'<2d$.  Then we have the following commutative diagram.
\begin{cd}
{} & \cHH^i(\fU,\sA_{2d}^*)\ar["\psi"]{dr}\ar["\sim"]{rr} & {} & \cHH^i(\fU,\sA_d^*)\\
\cHH^i(\fU,\sA_{2d'}^*)\ar{ur}\ar["\sim"]{rr} & {} & \cHH^i(\fU,\sA_{d'}^*)\ar{ur} & {}
\end{cd}
In particular, $\psi$ is injective and surjective, so an isomorphism, which implies the result.
\end{Proof}

We have now established the following sequence of morphisms, whose composition is the map $\varphi$ from Lemma \ref{sectionexists}.  We leverage Lemma \ref{cechOK} and the fact that colimits of abelian sheaves are exact.
\begin{eqnarray*}
\Pic(X)&\cong&\HH^1(X,\sO_X^*)\\
&\cong&\HH^1(X,\sO_X^{+*})\\
&\cong&\HH^1(X,\ilim\sA_d^*)\\
&\longto&\ilim\HH^1(X,\sA_d^*)\\
&\cong&\HH^1(X,\sA_d^*)\\
&\cong&\dlim\HH^1(X,\sA_d^*)\\
&\cong&\HH^1(X,\dlim\sA_d^*)\\
&\cong&\HH^1(X,\tilde\sO_X^*)\\
&\cong&\bZ[1/p].
\end{eqnarray*}
We make the necessary identifications to view the map $\varphi$ and its section $\sigma$ as maps between the following groups.
\begin{cd}
\HH^1(X,\sO_X^{+*}) \ar["\varphi"]{r} & \ilim\HH^1(X,\sA_d^*)\ar[dashed,bend right=30,swap,"\sigma"]{l}
\end{cd}
We can view the first group as isomorphism classes of invertible $\sO_X^+$ modules, and the second as inverse systems of isomorphism classes of invertible $\sA_d$ modules.  Under these identifications we have:
\begin{eqnarray*}
\varphi:&\sL\mapsto\{\sL/\varpi^d\sL\}\\
\sigma:&\{\sM_d\}\mapsto\ilim\sM_d.
\end{eqnarray*}
In particular, there is a natural map
\[\sL\longto\ilim\sL/\varpi^d\sL = \sigma\varphi\sL.\]
Locally, on an affinoid $\spa(R,R^+)$, we associate $\sL$ to an invertible $R$-module $M$.  Then this map becomes,
\[M\longto\ilim M/\varpi^d M\cong\hat M,\]
which is an isomorphism since $M$ is already complete.  We conclude that $\sigma$ is surjective, and therefore an isomorphism.  Putting all this together we may conclude that $\Pic X\cong\bZ[1/p]$, completing the proof of Theorem \ref{picpnperf}.

\subsection{Cohomology of Line Bundles}
In his Ph.D. thesis \cite{be18}, Bedi computed the cohomology of some of the twisting sheaves $\sO(d)$ on projectivoid space.  His proof was modeled on the computation for classical projective space in Ravi Vakil's algebraic geometry notes (\cite{FOAG} Chapter 18.2), but doesn't explicitly take into account the completions involved, and only includes the case for $n=2$.  Instead, we adapt the proof from EGA III \cite{EGA3} which relates \v{C}ech cohomology to the Koszul complex.  We also fill in the result, computing the cohomology of every line bundle in every degree.  Our rings are going to be topologically $\bZ[1/p]$ graded, and for a graded ring $A$ we will denote by $A_d$ the degree $d$ part.
\begin{Theorem}
Let $X = \bP^{n,\perf}$ be projectivoid space, and $\sO_X(d)\in\Pic X$ an arbitrary line bundle.  Then:
\begin{enumerate}
\item{
If $d\ge 0$,
\[\HH^0\left(X,\sO_X(d)\right) = K\left\la T_0^{1/p^\infty},\cdots,T_n^{1/p^\infty}\right\ra_d.\]
}
\item{
If $d<0$, then $\HH^n\left(X,\sO_X(d)\right)$ is the completion of the $K$ vector space generated by monomials of degree d, where the degree of each indeterminate is strictly negative, that is: 
\[\HH^n\left(X,\sO_X(d)\right) = \left\la T_0^{\alpha_0}\cdots T_n^{\alpha_n}\hspace{5pt}\big{|}\hspace{5pt} \alpha_i\in\bZ[1/p]_{<0}\text{ and } \sum\alpha_i=d\right\ra\widehat{ }\]
}
\item{
In all other cases,
\[\HH^r\left(X,\sO_X(d)\right) = 0.\]
}
\end{enumerate}
In particular:
\[ h^r\left(X,\sO_X(d)\right) =
\begin{cases}
1 & r = d = 0\\
\infty & r=0\text{ and }d>0\\
\infty & r=n\text{ and }d<0\\
0 & \text{ all other cases}
\end{cases}
.\]
\end{Theorem}
\begin{Proof}
We will leverage that colimits of abelian groups are exact, so that cohomology of finite complexes of abelian groups commute with arbitrary direct sums.  We will therefore study the \v{C}ech sequence associated to the sheaf
\[\sH = \bigoplus_{d\in\bZ[1/p]}\sO_X(d).\]
Let $A = K\left\la T_0^{1/p^\infty},\cdots,T_n^{1/p^\infty}\right\ra$, and $\fU = \{U_i\}$ be the standard cover of $X$.
\[\cech^r(\fU,\sH) = \prod_{0\le i_1<\cdots<i_r\le n}\widehat{A_{T_{i_1}\cdots T_{i_r}}}.\]
Since the differentials commute with degree and cohomology commutes with direct sums, we can conclude that
\[\cHH^i(\fU,\sO_X(d)) = \cHH^i(\fU,\sH)_d.\]
Furthermore, all finite intersections of the $U_i$ are affinoid, and vector bundles on affinoids are acyclic (\cite{ke17} Theorem 1.4.2).  Then the \v{C}ech to derived spectral spectral sequence degenerates to an isomorphism between \v{C}ech cohomology and sheaf cohomology, so we have reduced to computing these \v{C}ech cohomology groups.  We first consider the 0th cohomology. 
\[\cHH^0(\fU,\sH) = \bigcap_{i=0}^n\widehat{A_{T_i}} = A.\]
This proves the first statement of the theorem.  For the second, we consider the sequence
\[C^*(A):\hspace{10pt} 0\longto \prod A_{T_i} \longto \prod A_{T_iT_j}\longto\cdots\longto A_{T_0\cdots T_n}\longto0.\]
In each case, there is a countable fundamental system of neighborhoods of zero, given by $(\varpi^n)$, so that proceeding by induction from the left to the right, we see that completion on this sequence commutes with taking cohomology (see, for example, \cite{stacks-project} tag 0AMQ).   In particular,
\[\cHH^i(\fU,\sH) = \widehat{\HH^i\left(C^*(A)\right)}.\]
Let's analyze $H^n(C^*(A))$.  Notice that $A_{T_0\cdots T_n}$ is the $K$ vector space generated by monomials $T_0^{\alpha_0}\cdots T_n^{\alpha_n}$ for $\alpha\in\bZ[1/p]$.  The image of the $(n-1)$st differential is the $K$ vector space generated by monomials where at least one of the $\alpha_i\ge0$.  Therefore $H^n(C^*(A))$, which is the cokernel of this differential, is the $K$ vector space generated by monomials where each $\alpha_i<0$.  Taking completions proves the second statement of the theorem.

For the third statement of the theorem, the cases of $r<0$ and $r>n$ are trivial, so we assume $0<r<n$.  We will show that $H^r(C^*(A))=0$ since the completion of 0 is 0.  We point out that for any $f\in A$,
\[ A_{f} \cong \dlim\left(A\longtoo{\cdot f} A\longtoo{\cdot f}\cdots\right).\]
For all $s\ge0$, let $T^s = (T_0^s,\cdots,T_n^s)$.  Then $T^s$ is an $A$-regular sequence, and the associated Koszul complex $K^*(T^s)$ is a free resolution of $A/(T_0^s,\cdots,T_n^s)$.
\[K^*(T^s):\hspace{10pt} 0 \longto \Lambda^{n+1}A^{n+1}\longto\cdots\longto\Lambda^2 A^{n+1}\longto A^{n+1}\longtoo{(T_0^s,\cdots,T_n^s)} A\longto0.\]
In particular, the homology groups $\HH_i(K^*(T^s)) = 0$ for all $i>0$.  For each $s$ we can also look at the dual Koszul complex, and take the colimit as $s$ goes to infinity, with the above identification to the localization in mind.
\begin{cd}
{} & {} & \vdots\dar & \vdots\dar & \vdots\dar &{} & \vdots\dar & {}\\
C_s: & 0\rar & A\dar\ar["T^s"]{r} & A^{n+1}\ar["\cdot T"]{d}\ar["\cdot\wedge T^s"]{r} & \Lambda^2 A^{n+1}\ar["\cdot(T\wedge T)"]{d}\rar & \cdots\rar & \Lambda^{n+1}A^{n+1}\ar["\cdot T^{\wedge(n+1)}"]{d}\rar & 0\\
C_{s+1}: & 0\rar & A\dar\ar["T^{s+1}"]{r} & A^{n+1}\ar["\cdot\wedge T^{s+1}"]{r}\dar & \Lambda^2 A^{n+1}\rar\dar & \cdots\rar & \Lambda^{n+1} A^{n+1}\rar\dar & 0\\
{} & {} & \vdots\dar & \vdots\dar & \vdots\dar &{}&\vdots\dar&{}\\
\dlim C_s: & 0\rar &A\rar & \prod A_{T_i}\rar & \prod A_{T_iT_j}\rar & \cdots\rar & A_{T_0\cdots T_n}\rar & 0
\end{cd}
For $i>0$, the bottom row is $C^*(A)[1]$.  By the self-duality of the Koszul complex, we have that
\[\HH^i(C_s)\cong \HH_{n+1-i}(K^*(T^s)),\]
so that in particular, for $i<n+1$, we have $H^i(C_s) = 0$.  Since colimits of finite complexes commute with cohomology, we conclude that for $0<r<n$, 
\begin{eqnarray*}
\HH^r(C^*(A)) &\cong& \HH^{r+1}(\dlim C_s)\\
&\cong&\dlim \HH^{r+1}(C_s)\\
&=&\dlim 0\\
&=& 0.
\end{eqnarray*}
Taking completions proves the third statement of the theorem, and so we are done.
\end{Proof}
\begin{Remark}\label{integralcohomology}
An identical argument, but replacing $K$ with $K^\circ$, computes the \v Cech cohomology of the integral line bundles $\sO_X^+(d)$ for all $d>0$.  In particular we see that $\sO_X^+$ is \v Cech acyclic with respect to the standard cover, which completes the proof of Lemma \ref{cechacyclic}.

Since $\sO_X^+$ is only a priori \emph{almost} acyclic on affines (rather than acyclic on the nose like $\sO_X$), we cannot conclude from this that $\sO_X^+$ is acyclic on projectivoid space.  It is perhaps true that on the affinoids arising in the standard cover the sheaf of integral elements is indeed acyclic, and such a result would certainly simplify the arguments of Section \ref{computingpic}, but so far a proof has eluded the author.
\end{Remark}


\section{Maps to Projectivoid Space}\label{projectivoidgeometry}
Suppose $S$ is a scheme over $K$.  Then well known correspondence between maps from $S\to\bP^n$ and globally generated line bundles on $S$ together with a choice of $n+1$ generating global sections (see for example \cite{ha77} Theorem II.7.1).  In this section we will prove an analog of this correspondence for perfectoid spaces.

\begin{Definition}\label{newgroupoid}
To a perfectoid space $X$ over $K$, we associate a groupoid $\fL_n(X)$ whose objects consist of tuples $\left(\sL_i,s^{(i)}_j,\alpha_i\right)$ for $i\ge0$ and $j=0,\cdots,n$, where $\sL_i$ are line bundles on $X$, $s^{(i)}_0,\cdots,s^{(i)}_n\in\Gamma(X,\sL_i)$ are generating global sections, and $\alpha_i:\sL_{i+1}^{\otimes p}\longtoo{\sim}\sL_i$ are isomorphisms mapping $\left(s^{(i+1)}_j\right)^{\otimes p}\mapsto s_j^{(i)}$.  Morphisms are isomorphisms of line bundles which are compatible with the global sections and isomorphisms $\alpha_i$.

If $f:X\to Y$ is a $K$-morphism, we get a pullback functor $f^*:\fL_n(Y)\to\fL_n(X)$, so that $\fL_n$ is a category fibered in groupoids.
\end{Definition}
\begin{Remark}
We will show that $\fL_n$ is actually representable by a perfectoid space, and therefore a stack for the analytic topology.
\end{Remark}
\begin{Remark}
Note that if some $\alpha_i$ exists, it is unique.  Indeed, for each $i$ the global sections $s_j^{(i)}$ generate $\sL_i$, so that an isomorphism $\sL_{i+1}^{\otimes p}\longto\sL_i$ shows that the global sections  $\left(\s_{j}^{(i+1)}\right)^{\otimes p}$ generate $\sL_{i+1}^{\otimes p}$.  In particular, the isomorphism is completely determined by the images of these global sections.
\end{Remark}
\begin{Remark}
For each $i$, the data $\left(\sL_i,s^{(i)}_j\right)$ corresponds to a map to a projective space (as a rigid analytic variety), so that objects of the category $\fL_n(X)$ correspond to $p$th power root systems of maps to projective space. 
\end{Remark}
The main result of this section is that the category $\fL_n(X)$ parametrizes $K$-morphisms $X\longto\bP^{n,\perf}$.  In particular, viewing $\fL_n$ as a functor to sets we construct a natural isomorphism $\Hom(\cdot,\bP^{n,\perf})\cong\fL_n$ of functors from perfectoid spaces over $K$ to sets.  First we introduce a bit of notation.
\begin{Notation}
Denote by $m_i:\sO(1/p^{i+1})^{\otimes p}\longtoo{\sim}\sO(1/p^i)$ the isomorphism of line bundles on $\bP^{n,\perf}$ coming from multiplying factors together.
\end{Notation}
We now state the main theorem of this section (compare to \cite{ha77} Theorem II.7.1).
\begin{Theorem}\label{naturaliso}
The functor $\fL_n$ is represented by projectivoid space.

Explicitly, the natural transformation $\Hom(\cdot,\bP^{n,\perf})\to\fL_n$, which evaluated on $X$ takes $\phi:X\to\bP^{n,\perf}$ to the tuple $\left(\phi^*\sO(1/p^i),\phi^*T_j^{1/p^i},\phi^*m_i\right)\in\fL_n(X)$ is an isomorphism of functors.
\end{Theorem}
Since $\left\{T_j^{1/p^i}\right\}_{j=0}^n$ generates $\sO(1/p^i)$, we have that $\left\{\phi^*\left(T_j^{1/p^i}\right)\right\}_{j=0}^n$ generates $\phi^*\left(\sO(1/p^i)\right)$.  Furthermore, the standard isomorphisms $\sO(1/p^{i+1})^{\otimes p}\longtoo{\sim}\sO(1/p^i)$ coming from multiplying factors together send $(T_j^{1/p^{i+1}})^{\otimes p}$ to $T_j^{1/p^i}$, so pulling back these isomorphisms along $\phi$ gives us an element of $\fL_n(X)$.  We construct an inverse to this transformation in Section \ref{mainjuice8}, but first we will need a bit of setup.

\subsection{$\sL$-Distinguished Open Sets}
For this section we let $X$ be an adic space, $\sL$ a line bundle on $X$, and $s_1,\cdots,s_n$ global sections of $\sL$ which generate it at every point.  Let $D(s_i) = \{x\in X:s_i|_x$ generates $\sL_x\}$ be the \emph{doesn't vanish} set of the section $s_i$. Then the map $\sO_X\to\sL$ determined by $s_i$ is an isomoprhism on the stalks of every point of $D(s_i)$, and therefore restricts to an isomorphism on $D(s_i)$.  We suggestively denote the inverse by $s\mapsto s/s_i$.  Let's validate this notation with the following lemma.
\begin{Lemma}\label{relation}
On $D(s_i)\cap D(s_j)$, we have the following relation.
\[\frac{s_i}{s_j}\cdot\frac{s_j}{s_i} = 1.\]
\end{Lemma}
\begin{Proof}
We have two isomorphisms,
\begin{cd}
\Gamma(D(s_i)\cap D(s_j),\sO_X)\ar["s_{i}",yshift=3pt]{r}\ar[yshift = -3pt,swap, "s_j"]{r} & \Gamma(D(s_i)\cap D(s_j),\sL).
\end{cd}
Then we have
\[\frac{s_{i}}{s_{j}} = s_{j}^{-1}\circ s_{i}(1)\]
\[\frac{s_{j}}{s_{i}} = s_{i}^{-1}\circ s_{j}(1).\]
Since the maps $s_{i}^{-1}\circ s_{j}$ and $s_{j}^{-1}\circ s_{i}$ are inverses to each other, we win.
\end{Proof}
For every $x\in D(s_i)$, we can use the isomorphism $s_i^{-1}$ to get a valuation on $\Gamma(X,\sL)$.
\begin{cd}
\Gamma(X,\sL)\ar["res"]{r} &\Gamma(D(s_i),\sL)\ar["s_i^{-1}"]{r} & \Gamma(D(s_i),\sO_X)\ar["x"]{r} & \Gamma_x\cup\{0\}\\
s\ar[mapsto]{rrr} &{}&{}& {}|(s/s_i)(x)|
\end{cd}
With this in hand, we can define the following open subsets of $D(s_i)$ for each $i$.
\begin{Definition}
Let $X$ be a perfectoid space, $\sL$ a line bundle on $X$ and $s_1,\cdots,s_n$ generating global sections of $\sL$.  An open set of $X$ is called an \emph{$\sL$-distinguished open set} if it is of the form
\[X\left(\frac{s_1,\cdots,s_n}{s_i}\right) = \{x\in D(s_i):|(s_j/s_i)(x)|\le 1\text{ for all }j\}.\]
\end{Definition}
For the case of classical projective space, we can build a map to projective space along the \emph{doesn't vanish} sets of the given sections, and glue them together.  In the analytic topology these are not affinoid, so we must use these smaller $\sL$-distinguished open sets.  Let's prove these smaller open sets cover $X$.  Indeed, our notation suggests that one of $|(s_j/s_i)(x)|$ or $|(s_i/s_j)(x)|$ should be less than 1, let's check the details.
\begin{Lemma}
The $\sL$-distinguished open sets $X_i = X\left(\frac{s_1,\cdots,s_n}{s_i}\right)$ for $i=1,\cdots,n$ are open and cover $X$.
\end{Lemma}
\begin{Proof}
The openness of $X_i$ follows because it is in fact a \emph{rational} open in the adic space $D(s_i)$, which is open in $X$.  To show these cover $X$, fix some $x\in X$.  We already know the $D(s_i)$ cover $X$, because the $s_i$ generate $\sL$.  Therefore $I = \left\{i\in\{1,\cdots,n\}:x\in D(s_i)\right\}$ is nonempty and finite. Order the elements of $I$ via $i\le j$ if $|(s_i/s_j)(x)|\le 1$.  Notice that for any $i,j\in I$, we have $i\le j$ or $j\le i$.  Indeed, applying Lemma \ref{relation} together with the multiplicativity of the valuation given by $x$, we have either $|(s_i/s_j)(x)|\le 1$ or $|(s_j/s_i)(x)|\le 1$.  Also if $i\le j$ and $j\le k$ then
\[|(s_i/s_k)(x)| = |(s_i/s_j)(x)|\cdot|(s_j/s_k)(x)|\le 1,\]
so that $i\le k$.  Finally, notice that if $i\le j$ and $j\le i$ then we have $|(s_i/s_j)(x)| = |(s_j/s_i)(x)| = 1$.  Therefore we can choose (not necessarily uniquely) some $r$ which is maximal with respect to this ordering.  Then $|(s_i/s_r)(x)|\le 1$ for all $i\in I$.  For all other $i$, we have $x\notin D(s_i)$ so that $|(s_i/s_r)(x)| = 0 \le 1$.  Therefore $x\in X_r$, completing the proof.
\end{Proof}
\begin{Example}
The standard cover of $\bP^{n,\perf}$ by perfectoid unit disks consists of the $\sO(1)$-distinguished open sets $\bP^{n,\perf}\left(\frac{T_0,\cdots,T_n}{T_i}\right)$.
\end{Example}
The following lemma which implies that if $\sL^{\otimes p}\cong \sM$ and $s$ and $t$ are global sections of $\sL$ and $\sM$ respectively, with $s^{\otimes p} = t$, then $D(s) = D(t)$.  In particular, using the notation of Definition \ref{newgroupoid}, this implies that if the $s^{(i)}_j$ generate $\sL_i$, then the $s^{(i+1)}_j$ generate $\sL_{i+1}$.
\begin{Lemma}\label{generators}
Let $(R,\fm)$ be a local ring, and $M,N$ invertible $R$-modules such that $M^{\otimes r}\cong N$ for a positive integer $r$.  Let $f\in M$, and $g\in N$ such that under this identification $f^{\otimes r} = g$.  Then if $g$ generates $N$, $f$ generates $M$.
\end{Lemma}
\begin{Proof}
We show the contrapositive.  If $f$ does not generate $M$, then by Nakayama's lemma, $f\in\fm M$.  Thus $f=a\cdot s$ for some $a\in\fm$ and $s\in M$.  But then under the appropriate identification,
\[g = f^{\otimes r} = (a\cdot s)^{\otimes r} = a^r\cdot s^{\otimes r}\in\fm^r N\subseteq\fm N.\]
Therefore $g$ cannot generate $N$.
\end{Proof}

\subsection{Construction of the Projectivoid Morphism}\label{mainjuice8}
We can now finish the proof of Theorem \ref{naturaliso} by constructing an inverse to the natural tranformation from the theorem.  The result follows from the following proposition.
\begin{Proposition}\label{buildingthemap}
Let $X$ be a perfectoid space over $K$ and $\left(\sL_i,s_j^{(i)},\alpha_i\right)\in\fL_n(X)$.  There is a unique $K$-morphism $\phi:X\longto\bP^{n,\perf}$ such that
\[\left(\phi^*\sO(1/p^i),\phi^*T_j^{1/p^i},\phi^*m_i\right)\cong\left(\sL_i,s_j^{(i)},\alpha_i\right).\]
\end{Proposition}
\begin{Proof}
Let $X_j = X\left(\frac{s_0^{(0)},\cdots,s_n^{(0)}}{s^{(0)}_j}\right)$ be the cover of $X$ by $\sL_0$-distinguished opens.  Let $U_j = \bP^{n,\perf}\left(\frac{T_0,\cdots,T_n}{T_j}\right)\subseteq\bP^{n,perf}$ be the standard cover by affinoids.  The $U_i$ are isomorphic to the perfectoid unit polydisk, and are naturally identified with
\[\spa\left(K\left\la\left(\frac{T_0}{T_j}\right)^{1/p^\infty},\cdots,\left(\frac{T_n}{T_j}\right)^{1/p^{\infty}}\right\ra,K^\circ\left\la\left(\frac{T_0}{T_j}\right)^{1/p^\infty},\cdots,\left(\frac{T_n}{T_j}\right)^{1/p^{\infty}}\right\ra\right).\]
We build $\phi$ locally from maps $\phi_j:X_j\to U_j$.  Since $U_j$ is affinoid, it is equivalent to build a map of Huber pairs,
\[\left(K\left\la\left(\frac{T_0}{T_j}\right)^{1/p^\infty},\cdots,\left(\frac{T_n}{T_j}\right)^{1/p^{\infty}}\right\ra, K^\circ\left\la\left(\frac{T_0}{T_j}\right)^{1/p^\infty},\cdots,\left(\frac{T_n}{T_j}\right)^{1/p^{\infty}}\right\ra\right)\longtoo{\gamma_j}\left(\sO_X(X_j),\sO_X^+(X_j)\right).\]
That is, a ring map
\[K\left\la\left(\frac{T_0}{T_j}\right)^{1/p^\infty},\cdots,\left(\frac{T_n}{T_j}\right)^{1/p^{\infty}}\right\ra\longtoo{\gamma_j}\Gamma(X_j,\sO_X),\]
satisfying
\[\gamma_j\left(K^\circ\left\la\left(\frac{T_0}{T_j}\right)^{1/p^\infty},\cdots,\left(\frac{T_n}{T_j}\right)^{1/p^{\infty}}\right\ra\right)\subseteq\Gamma\left(X_j,\sO_X^+\right).\]
We define $\gamma_j$ on generators by the rule
\[\gamma_j\left(\left(\frac{T_r}{T_j}\right)^{1/p^i}\right) = \frac{s^{(i)}_{r}}{s_{j}^{(i)}}.\]
To make sure this is a ring homomorphism, we must check that
\[\left(\frac{s^{(i+1)}_r}{s^{(i+1)}_{j}}\right)^p = \frac{s^{(i)}_{r}}{s^{(i)}_{j}}.\]
First notice that, under the identification $\alpha_i:\sL_{i+1}^{\otimes p}\cong \sL_i$, the following diagram commutes (keeping in mind that the horizontal maps are not homomorphisms).
\begin{cd}
\sO_X \ar["x\mapsto x^p"]{r}\ar["s_{j}^{(i+1)}"]{d} & \sO_X \ar["s^{(i)}_{j}"]{d}\\
\sL_{i+1}\ar["s\mapsto s^{\otimes p}"]{r} & \sL_i.
\end{cd}

Indeed, the commutativity of this diagram follows directly from the multilinearity of tensor product together with the identification $\left(s_{j}^{(i+1)}\right)^{\otimes p} = s_{j}^{(i)}$.  Chasing this diagram, we see that
\begin{eqnarray*}
\left(\frac{s_{r}^{(i+1)}}{s_{j}^{(i+1)}}\right)^p &=& \left(\left(s_{j}^{(i+1)}\right)^{-1}\left(s_{r}^{(i+1)}\right)\right)^p\\
&=&\left(s_{j}^{(i)}\right)^{-1}\left(s_{r}^{(i)}\right)\\
&=&\frac{s_{r}^{(i)}}{s_{j}^{(i)}},
\end{eqnarray*}
as desired.  Therefore $\gamma_j$ is a homomorphism.  Finally, the definition of $X_j$ implies that for all $x\in X_j$,
\[\left|\gamma_j\left(\frac{T_i}{T_j}\right)(x)\right| =\left|\frac{s_{i}^{(0)}}{s_{k}^{(0)}}(x)\right| \le 1,\]
so that 
\[\gamma_j\left(\frac{T_i}{T_j}\right)\in\Gamma\left(X_j,\sO_X^+\right).\]
The multiplicativity of the valuation associated to $x$ shows the same holds for all $p$th power roots so that 
\[\gamma_j\left(K^\circ\left\la\left(\frac{T_0}{T_j}\right)^{1/p^\infty},\cdots,\left(\frac{T_n}{T_j}\right)^{1/p^{\infty}}\right\ra\right)\subseteq\Gamma\left(X_j,\sO_X^+\right).\]
Therefore we get a morphism $\phi_j:X_j\to U_j\subseteq\bP^{n,\perf}$, for each $j$.  Notice also that $s_{r}^{(i)}/s_{j}^{(i)}$ is a $p^i$-th root of $s_{r}^{(0)}/s_{j}^{(0)}$ as desired.

Finally, we check that these morphisms glue to a map $\phi:X\to\bP^{n,\perf}$.  This amounts to showing that the restrictions of $\gamma_j$ and $\gamma_k$ are equal as maps from $\Gamma(U_j\cap U_k,\sO_{\bP^{n,\perf}})\longto\Gamma(X_j\cap X_k,\sO_X)$.  That is, that
\[\gamma_j\left(\left(\frac{T_k}{T_j}\right)^{1/p^i}\right) = \gamma_k\left(\left(\frac{T_j}{T_k}\right)^{1/p^i}\right)^{-1}.\]
With our notation, this boils down to
\[\frac{s_{k}^{(i)}}{s_{j}^{(i)}}\cdot\frac{s_{j}^{(i)}}{s_{k}^{(i)}} = 1.\]
But this is just Lemma \ref{relation}.

The rest is immediate from the construction.  Since $\sO_{\bP^{n\perf}}(d)$ is generated by the monomials of degree $d$, the construction shows that $\phi^*\sO(1/p^i) = \sL_i$, and that $\phi^*\left(T_j^{1/p^i}\right) = s_{j}^{(i)}$.  Furthermore, any map $\psi:X\to\bP^{n,\perf}$ with these properties is by definition given affinoid locally on the standard cover of the target by the ring map $(T_j/T_k)^{1/p^i}\mapsto s^{(i)}_j/s^{(i)}_k$.  That is, $\psi|_{X_i}=\phi_i$ so that $\psi=\phi$.
\end{Proof}
\subsection{The Positive Characteristic Case}
If $X$ is a perfectoid space of characteristic $p$, then the Frobenius morphism $Frob:\bG_m\to\bG_m,x\mapsto x^p$ is an isomorphism.  Therefore the $p$th power map on $\Pic X$ is an isomorphism as well, since it is $H^1(X,Frob)$.  This means that given $\left(\sL_i,s^{(i)}_j,\alpha_i\right)\in\fL_n(X)$, the $\sL_i$ for $i>0$ are uniquely determined by $\sL_0$.  Similarly, since $X$ is perfect, the map $\gamma_i$ constructed in the proof of Proposition \ref{mainjuice8} is completely determined by where $T_r/T_i$ goes for each $r\not=j$, because the $p$th roots of the image are unique.  We summarize this in the following corollary.
\begin{Corollary}\label{FOPCP}
If $X$ is a perfectoid space over $K$ of characteristic $p$, a map $X\to\bP^{n,\perf}$ is equivalent to a line bundle $\sL$ on $X$ and global sections $s_0,\cdots,s_n$ that generate $\sL$, or equivalently, to a map to classical projective space $\bP^n$ (as a rigid analytic variety).
\end{Corollary}
We can now leverage the tilting equivalence to say that maps to $X\to\bP^{n,\perf}$ in any characteristic are governed by a single line bundle on $X^\tilt$.  Indeed, by the tilting equivalence we have that $\Hom(X,\bP^{n,\perf}_K) = \Hom(X^\tilt,\bP^{n,\perf}_{K^\tilt})$. This implies the following corollary to Theorem \ref{naturaliso}.
\begin{Corollary}\label{tiltmaps}
If $X$ is a perfectoid space over $K$ of any characteristic, a map $X\to\bP^{n,\perf}_K$ is equivalent to a single line bundle $\sL$ on $X^\tilt$ together $n+1$ global sections generating $\sL$.
\end{Corollary}
Using this corollary as an intermediary, we get a natural and geometric correspondence between certain inverse systems of line bundles on $X$ and single line bundles on $X^\tilt$.
\begin{Corollary}\label{bundlesandtilting}
An element of $\fL_n(X)$ is equivalent to the a line bundle $\sL\in\Pic X^\tilt$ together with $n+1$ generating global sections.
\end{Corollary}
This will be a useful tool in understanding the relationship between $\Pic X$ and $\Pic X^\tilt$.



\section{Untilting Line Bundles}\label{untiltingLB}
Recall that one of our motivations was to understand the behavior of vector bundles under the tilting equivalence.  In this section, we use the tools of projectivoid geometry developed in Section \ref{projectivoidgeometry} to compare the Picard groups of a perfectoid space $X$ and its tilt $X^\tilt$.  Indeed, the theory of maps to projectivoid space allows us to pass between line bundles on $X$ and $X^\tilt$ by choosing (compatible) generating sections, constructing the associated map to projectivoid space, and then using the tilting equivalence to pass across characteristics.  We remark that the theory of \emph{pro-\'etale cohomology} on perfectoid spaces allows us to make this comparison cohomologically, but the geometric theory we developed in the previous section gives us a firm geometric grasp.

\subsection{Cohomological Untilting}
In \cite{bhash15}, Bhatt and Scholze introduce the \emph{pro-\'etale} site for schemes and perfectoid spaces.  We review the definition here.
\begin{Definition}
A map $f:Y=\spa(S,S^+)\to X=\spa(R,R^+)$ of affinoid perfectoid spaces is called \emph{affinoid pro-\'etale} if it can be written as a cofiltered limit of \'etale maps $Y_i = \spa(S_i,S_i^+)\to X$ of affinoid perfectoid spaces.  More generally, a map $f:Y\to X$ of perfectoid spaces is \emph{pro-\'etale} if is locally on the source and target affinoid pro-\'etale.

The (small) pro-\'etale site of $X$ is the Grothendieck topology  on the category of perfectoid spaces $f:Y\to X$ pro-\'etale over $X$ on which a collection $\{f_i:Y_i\to X\}_{i\in I}$ is a covering if for each quasicompact open $U\subseteq X$ there exists a finite subset $J\subseteq I$ and quasicompact open subsets $V_i\subseteq Y_i$ for $i\in J$ such that $U = \cup_{i\in J}f_i(V_i)$.

If $\sF$ is a pro-\'etale sheaf on $X$ (that is a sheaf on the pro-\'etale site of $X$), the pro-\'etale cohomology groups $\HH^i(X_\proet,\sF)$ are the derived functor sheaf cohomology groups on the pro-\'etale site.
\end{Definition}
Let $X$ be a perfectoid space over $K$.  The pro-\'etale sheaf $\bG_{m,X}$ maps $U\mapsto\Gamma(U,\sO_U)^*$.  We have the following theorem.
\begin{Proposition}\label{proetalepicard}
$\HH^1(X_\proet,\bG_m)\cong\Pic X$.
\end{Proposition}
\begin{Proof}
For any site $S$, the cohomology group $\HH^1(X_S,\bG_m)$ parametrizes isomorphism classes of line bundles on $X$ with respect to the topology of $S$.  Due to \cite{KL16} Theorem 3.5.8, vector bundles (of any rank) on a perfectoid space with respect to the pro-\'etale, \'etale, and analytic topologies coincide.
\end{Proof}
We use the equivalence of the pro-\'etale topologies of $X$ and $X^\tilt$ to construct the tilt of $\bG_m$ as a pro-\'etale sheaf on $X$:
\[\bG_{m,X}^\tilt:U\mapsto (\Gamma(U,\sO_U)^\tilt)^* = \Gamma(U^\tilt,\sO_{U^\tilt})^* = \Gamma(U^\tilt,\bG_{m,X^\tilt}).\]
The equivalence of the \'etale topologies on $X$ and $X^\tilt$ show that $\bG_{m,X}^\tilt$ is indeed a sheaf.  Better yet, the effectiveness of \v{C}ech cohomology on the pro-\'etale site shows that
\[\HH^i(X_\proet,\bG_{m,X}^\tilt)\cong \HH^i(X^\tilt_\proet,\bG_{m,X^\tilt}).\]
In particular, $\HH^1(X_\proet,\bG_{m,X}^\tilt)\cong\Pic X^\tilt$.
Now consider the Kummer sequence for various powers of $p$.
\[0\longto\bmu_{p^n}\longto\bG_{m,X}\longto\bG_{m,X}\longto0.\]
This is an exact sequence of sheaves on the pro-\'etale site of $X$.  Indeed, this can be checked on the stalks, which on the pro-\'etale site are strictly Henselian local rings.  Therefore we can form an inverse system of exact sequences:
\begin{cd}
0\rar & \bmu_p\rar & \bG_{m,X} \rar & \bG_{m,X}\rar & 0\\
{} & \vdots \ar{u} & \vdots\ar{u} & \vdots\ar[equals]{u} & {}\\
0\rar & \bmu_{p^n}\rar\ar{u} & \bG_{m,X}\rar\ar{u} & \bG_{m,X}\rar\ar[equals]{u} & 0\\
0\rar & \bmu_{p^{n+1}}\rar\ar{u} & \bG_{m,X}\rar\ar{u} & \bG_{m,X}\rar\ar[equals]{u} & 0\\
{} & \vdots \ar{u} & \vdots\ar{u} & \vdots\ar[equals]{u} & {}\\
\end{cd}
The vertical maps on the the left and middle sides are $x\mapsto x^p$.  Taking this limit gives the following sequence.
\[0\longto\bZ_p(1)\longto\bG_{m,X}^\tilt\longtoo{\sharp}\bG_{m,X}.\]
The middle term is $\bG_{m,X}^\tilt$ essentially by definition.  Indeed, the construction the tilt of a perfectoid algebra $R$ (\cite{sh12} Lemma 3.4) induces an map of multiplicative monoids:
\[R^\tilt\cong\ilim_{x\mapsto x^p} R,\]
which restricts to the desired isomorphism on unit groups.  Finally, exactness on the right can be checked explicitly in the pro-\'etale topology.  Indeed, adjoining a $p$th power root is \'etale cover so that adjoining all the missing $p$-power roots gives a pro-\'etale cover on which $\sharp$ is surjective.  Therefore we have a short exact sequence of pro-\'etale sheaves:
\[0\longto\bZ_p(1)\longto\bG_{m,X}^\tilt\longtoo{\sharp}\bG_{m,X}\longto0.\]
\begin{Remark}\label{teich}
If $R$ is a perfectoid algebra we always get a map of monoids $\sharp:R^\tilt\to R$ given by projection onto the first coordinate.  Although it is not a ring homomorphism unless $R$ already had characteristic $p$, its restriction to unit groups $(R^\tilt)^*\to R^*$ is a group homomorphism.  This construction is another way of building the map $\sharp:\bG_{m,X}^\tilt\to\bG_{m,X}$.  The advantage of the above construction is that it explicitly exhibits the Tate module $\bZ_p(1)$ as the kernel.
\end{Remark}
Taking long exact sequences in cohomology gives us the following diagram, where the rows are exact.
\begin{cd}
{} & \vdots & \vdots & \vdots & \vdots & {}\\
\cdots\rar & \HH^1(X_\proet,\bmu_{p^n})\rar\ar{u} & \Pic X\rar\ar{u} & \Pic X\rar\ar[equals]{u} & \HH^2(X_\proet,\bmu_{p^n})\rar\ar{u}&\cdots\\
\cdots\rar & \HH^1(X_\proet,\bmu_{p^{n+1}})\rar\ar{u} & \Pic X\rar\ar["\sL\mapsto\sL^{\otimes p}"]{u} & \Pic X\rar\ar[equals]{u} & \HH^2(X_\proet,\bmu_{p^{n+1}})\rar\ar{u}&\cdots\\
{} & \vdots\ar{u} & \vdots\ar{u} & \vdots\ar[equals]{u} & \vdots\ar{u} & {}\\
\cdots\rar & \HH^1(X_\proet,\bZ_p(1))\rar\ar{u} & \Pic X^\tilt\ar["\theta_0"]{r}\ar{u}\ar[bend right=30,swap,"\theta_{n+1}"]{uu}\ar[bend right=60,swap,"\theta_n"]{uuu} & \Pic X\rar\ar[equals]{u} & \HH^2(X_\proet,\bZ_p(1))\rar\ar{u}&\cdots\\
\end{cd}
Taking the inverse limit of the $\theta_n$ gives a homomorphism of groups,
\begin{equation}\label{pi}
\theta:\Pic X^\tilt\longto\ilim_{\sL\mapsto\sL^{\otimes p}}\Pic X,
\end{equation}
and $\theta_0$ is this map composed with the projection onto the first coordinate.
\begin{Remark}
In Corollary \ref{bundlesandtilting} we established that inverse systems of $p$th roots of line bundles (with generating sections) on $X$ correspond to individual line bundles (with generating sections) on $X^\tilt$.  This seems to suggest that $\theta$ could be an isomorphism in cases where we have nice maps to projective space.
\end{Remark}
\subsection{Untilting Via Maps to Projectivoid Space}
We now give a geometric understanding of $\theta$ and $\theta_0$ in terms of maps to projectivoid space, and in doing so study whether whether the following correspondence holds for a perfectoid space $X$.
\[
\left\{\begin{tabular}{p{3cm}}Globally generated $\sL\in\Pic X^\tilt$\end{tabular}\right\}\leftrightarrow\left\{\begin{tabular}{p{5cm}}Systems of globally generated line bundles $(\sL_0,\sL_1,\cdots)$ on $X$ such that $\sL_{i+1}^{\otimes p}\cong\sL_i$.\end{tabular}\right\}
\]
We begin by constructing a map in the righthand direction.  Given a globally generated invertible sheaf $\sL\in\Pic X^\tilt$, choose $n$ sections which generate $\sL$.  Associated to this data there is a unique morphism $\phi^\tilt:X^\tilt\longto\bP^{n,\perf}_{K^\tilt}$, which is the tilt of a unique morphism $\phi:X\longto\bP^{n,\perf}_K$.  Let $\sL_i=\phi^*(\sO(1/p^i))$.  This gives a system of $(\sL_0,\sL_1,\cdots)$ on the right hand side.  As a first step we show that the sheaves $\sL_i$ do not depend on the choices of global sections of $\sL$.
\begin{Proposition}\label{geometrywelldefined}
The construction in the previous paragraph is well defined, and $(\sL_0,\sL_1,\cdots) = \theta(\sL)$ where $\theta$ is the cohomological map defined in Equation \ref{pi} above.
\end{Proposition}
\begin{Proof}
$\phi^*$ can be constructed cohomologically by applying cohomology to the unit of the adjunction, $u:\bG_{m,\bP^{n,\perf}_K}\to\phi_*\bG_{m,X}$ and composing with the natural map $\HH^1(\bP^{n,\perf},\phi_*\bG_{m,X})\to\HH^1(X,\bG_{m,X})$.  Pulling $u$ back along the $p$th power map gives $\phi^{\tilt*}$ the same way.  Since the $p$th power map commutes with pullback, we get the following commutative diagram.
\begin{cd}
\Pic\bP^{n,\perf}_{K^\tilt}\ar["\phi^{\tilt*}"]{r}\ar["\theta_{\bP^{n,\perf}}"]{d} & \Pic X^\tilt\ar["\theta_X"]{d}\\
\ilim\Pic\bP^{n,\perf}_K\ar["\phi^*"]{r} & \ilim\Pic X.
\end{cd}
Since $\sL=\phi^{\tilt*}\sO_{\bP^{n,\perf}_{K^\tilt}}(1)$ and $\sL_i = \phi^*\sO_{\bP^{n,\perf}_K}(1/p^i)$, we have reduced to proving the proposition for $\bP^{n,\perf}_K$.  Explicitly, we must show
\[\theta_{\bP^{n,\perf}}(\sO(1)) = (\sO(1),\sO(1/p),\sO(1/p^2),\cdots).\]
Since $\Pic\bP^{n,\perf} = \bZ[1/p]$, and is therefore uniquely $p$-divisible, it is enough to show that
\[\theta_{0,\bP^{n,\perf}}\sO(1) = \sO(1).\]
Now $\theta_0$ is the cohomological map associated to the Teichmuller map $\sharp:\bG_m^\tilt\longto\bG_m$.  Scholze showed in \cite{sh12} Proposition 5.20 that the Teichmuller map on the perfectoid Tate algebra maps $T_i\mapsto T_i$.  View $\theta_0$ as a map on \v{C}ech cohomology with respect to the standard affine covers, and view $\HH^1(\bP^{n,\perf},\bG_m)$ as descent data for building a line bundle (and similarly for the tilt).  Then we see that $\sharp$ sends descent data for $\sO(1)$ (which is monomials of degree one), to monomials of degree one, which build $\sO(1)$ on $\bP^{n,\perf}_K$.

Therefore untilting line bundles via maps to projectivoid is a well defined process, as it agrees with the cohomological method which does not depend on the choice of sections.
\end{Proof}

We split the study of whether this bijects into the two cases of injectivity and surjectivity.  We first analyze each case and see where the difficulties may lie.
\begin{itemize}
\item{
\textbf{Injectivity}: Suppose $\sL,\sM\in\Pic(X^\tilt)$ are globally generated.  Suppose choosing sections and untilting the associated maps to projectivoid space gives us maps $\phi:X\to\bP^{n,\perf}_K$ and $\psi:X\to\bP^{r,\perf}_K$.  If $\phi^*(\sO(1/p^i))\cong\psi^*(\sO(1/p^i)=:\sL_i$ for all $i$, can we conclude that $\sL\cong\sM$?  We can attack this using the methods of Section \ref{projectivoidgeometry} by considering the tuples $(\sL_i,s^{(i)}_j,\alpha_i)\in\fL_n(X)$ and $(\sL_i,t^{(i)}_j,\beta_i)\in\fL_r(X)$ associated to $\phi$ and $\psi$ respectively.  If the $\alpha_i$ and $\beta_i$ agree, we can consider $\left(\sL_i,\left\{s^{(i)}_j,t^{(i)}_k\right\},\alpha_i\right)\in\sL_{n+r+1}(X)$ and consider how the associated map $X\to\bP^{n+r+1,\perf}_K$ tilts.  In this section we settle the case where $\alpha_i=\beta_i$, and solve the general case in Section \ref{injectiveTheta}.
}
\item{
\textbf{Surjectivity}: Suppose $(\sL_0,\sL_1,\cdots)$ are globally generated with $\sL_{i+1}^{\otimes p}\cong\sL_i$, and there are global sections $s_{j}^{(i)}$ generating $\sL_i$ such that $\left(s_{j}^{(i)}\right)^{\otimes p} = s_{j}^{(i)}$.  Then passing through the maps to projective space we get $\sL\in\Pic X^\tilt$ which maps to $(\sL_0,\sL_1,\cdots)$ under $\theta$.  But, can we always find sections $s_{j}^{(i)}$ and isomorphisms such that $\left(s_{j}^{(i+1)}\right)^{\otimes p} = s_{j}^{(i)}$?  Restated, are there generating global sections of $\sL_0$ all of whose $p$th power roots exist?  If so, our correspondence surjects.
}
\end{itemize}
In the rest of this section we settle injectivity in the case where the isomorphisms $\sL_{i+1}^{\otimes p}\cong\sL_{i}$ agree for the two sets of sections.
\begin{Proposition}\label{matchingisos}
Let $X$ be a perfectoid space over $K$.  Suppose $\left(\sL_i,s^{(i)}_j,\alpha_i\right)\in\fL_n(X)$ and $\left(\sL_i,t^{(i)}_j,\alpha_i\right)\in\fL_r(X)$ correspond to maps $\phi:X\to\bP^{n,\perf}_K$ and $\psi:X\to\bP^{r,\perf}_K$ respectively.  Then
\[\phi^{\tilt*}\sO_{\bP^{n,\perf}_{K^\tilt}}(1)\cong\psi^{\tilt*}\sO_{\bP^{r,\perf}_{K^\tilt}}(1).\]
\end{Proposition}

Fix $\left(\sL_i,s_j^{(i)},\alpha_i\right)$ corresponding to a map $\phi:X\to\bP^{n,\perf}_{K}$.  As a first step, we show that we can add one global section to each $\sL_i$ that are compatible with the $\alpha_i$ without changing the line bundle we get over $X^\tilt$.  Suppose $t_i\in\Gamma(X,\sL_i)$ is a global section such that $\alpha_i\left(t_{i+1}^{\otimes p}\right) = t_i$.  For every $\lambda = (\lambda_0,\lambda_1,\cdots)\in\ilim K^* = K^{\tilt*}$, we let $\psi_\lambda:X\longto\bP^{n+1,\perf}_K$ be the projectivoid map corresponding to adding $t_i$, that is,  corresponding to $\left(\sL_i,\{s_j^{(i)},t_i\},\alpha_i\right)$.   We hope to fit $\phi$ and $\psi_\lambda$ in a commutative diagram.  To do so we must develop an analog of rational maps in this analytic context.

If we want to define a map $\bP^{n+1,\perf}\longto\bP^{n,\perf}$ given by $\left(\sO(1/p^i),\left\{T_0^{1/p^i},\cdots,T_n^{1/p^i}\right\},m_i\right)$ we would notice that this isn't defined wherever $|T_i/T_{n+1}|>1$.  In particular, it is only defined on the open set:
\[U=\bigcup_{j\not=n+1}\bP^{n+1,\perf}_K\left(\frac{T_0,\cdots,T_{n+1}}{T_j}\right).\]
This is the projectivoid analog of projecting away from the hyperplane where $T_{n+1}$ vanishes, (here we are projecting away from a perfectoid disk at the `north pole').  Unfortunately, the image of $\psi_\lambda$ does not lie in $U$, because there may be points $x$ where $\left|(s_{j}^{(0)}/t_0)(x)\right|>1$ for all $i$, so that $|(T_i/T_{n+1})(\psi_\lambda(x))|>1$.  But, restricted to the open set
\[V_\lambda = \bigcup_j X\left(\frac{s_{0}^{(0)},\cdots,s_{n}^{(0)},t_0}{s_{j}^{(0)}}\right),\]
the image of $\psi_\lambda$ does lie in $U$.  Thus we have the following commutative diagram for every $\lambda$.
\begin{cd}
{} & {} & \bP^{n,\perf}_K\\
X\ar["\phi"]{urr}\ar[swap,"\psi_\lambda"]{drr}\ar[hookleftarrow]{r} & V_\lambda\ar["\psi_\lambda"]{r} & U\ar["\pi"]{u}\ar[hookrightarrow]{d}\\
{} & {} & \bP^{n+1,\perf}_K
\end{cd}
\begin{Lemma}\label{covering}
The sets $V_{(\varpi^{\tilt})^r}$ form an open cover of $X$.  As a consequence the sets $V_{(\varpi^\tilt)^r}^\tilt$ cover $X^\tilt$.
\end{Lemma}
\begin{Proof}
Notice $(\varpi^{\tilt})^r = (\varpi^r,\varpi^{r/p},\cdots)$.  Fix $x\in X$.  There is some $j$ such that $x\in X\left(\frac{s_{0}^{(0)},\cdots,s_{n}^{(0)}}{s_{j}^{(0)}}\right)$.  Furthermore, since $\varpi$ is topologically nilpotent, there is some $r$ such that
\[\left|(\varpi^r t_0/s_{j}^{(0)})(x)\right| = |\varpi^r|\cdot\left|(t_0/s_{j}^{(0)})(x)\right|<1\]
proving the first statement.  The second is an immediate consequence of the tilting equivalence.
\end{Proof}

\begin{Lemma}\label{rationalpullbacks}
For any $\lambda\in K^{\tilt *}$,
\[\left(\phi^{\tilt*}\sO_{\bP^{n,\perf}_{K^\tilt}}(1)\right)|_{V^\tilt_{\lambda}}\cong\psi_\lambda^{\tilt*}\left(\sO_{\bP^{n+1,\perf}_{K^\tilt}}(1)|_{U^\tilt}\right)\cong\left(\psi_\lambda^{\tilt*}\sO_{\bP^{n+1,\perf}_{K^\tilt}}(1)\right)|_{V_{\lambda}^\tilt}\]
\end{Lemma}
\begin{Proof}
This follows from the commutativity of the tilt of the diagram above, reproduced below, together with the fact that $\pi^\tilt$ is given by the line bundle $\sO_{\bP^{n+1,\perf}_{K^\tilt}}(1)|_U$ together with the sections $T_0,\cdots,T_n$.
\begin{cd}
{} & {} & \bP^{n,\perf}_{K^\tilt}\\
X^\tilt\ar["\phi^\tilt"]{urr}\ar[swap,"\psi^\tilt_\lambda"]{drr}\ar[hookleftarrow]{r}& V^\tilt_\lambda\ar["\psi^\tilt_\lambda"]{r} & U^\tilt\ar["\pi^\tilt"]{u}\ar[hookrightarrow]{d}\\
{} & {} & \bP^{n+1,\perf}_{K^\tilt}
\end{cd}
\end{Proof}
\begin{Lemma}\label{differentunits}
Fix any $\lambda,\xi\in\ilim K^* = K^{\tilt*}$.  Then
\[\psi^{\tilt*}_\lambda\sO_{\bP^{n+1,\perf}_{K^\tilt}}(1)\cong\psi^{\tilt*}_\xi\sO_{\bP^{n+1,\perf}_{K^\tilt}}(1).\]
\end{Lemma}
\begin{Proof}
Let $\tau:\bP^{n+1,\perf}_K\to\bP^{n+1,\perf}_K$ be the map determined by $\left(\sO(1/p^i),\left\{T_0^{1/p^i},\cdots,T_n^{1/p^i},\frac{\gamma_i}{\xi_i}T_{n+1}^{1/p^i}\right\},m_i\right)$.  Then $\tau$ is an isomorphism, and $\tau^\tilt$ is the map determined by $\sO(1)$ and $T_0,\cdots,T_n,\frac{\gamma}{\xi}T_{n+1}$.  We have the following two commutative diagrams, the right hand diagram being the tilt of the left.
\begin{cd}
{}&\bP^{n+1,\perf}_K \ar["\tau"]{dd}&{} &{} &\bP^{n+1,\perf}_{K^\tilt}\ar["\tau^\tilt"]{dd}\\
X\ar["\psi_\lambda"]{ur}\ar[swap,"\psi_\xi"]{dr} & {} & {} &X^\tilt\ar["\psi_\lambda^\tilt"]{ur}\ar[swap,"\psi_\xi^\tilt"]{dr} &{}\\
{} &\bP^{n+1,\perf}_K&{}&{}&\bP^{n+1,\perf}_{K^\tilt}
\end{cd}
Since $\tau^{\tilt*}\sO(1)=\sO(1)$, we are done.
\end{Proof}
Putting these three lemmas together, we conclude that
\[\phi^{\tilt*}\sO_{\bP^{n,\perf}_{K^\tilt}}(1)\cong\psi_1^{\tilt*}\sO_{\bP^{n+1,\perf}_{K^\tilt}}(1).\]
Indeed, the pullback of $\sO(1)$ along $\psi_1^\tilt$ agrees with the pullback along $\psi_{(\varpi^\tilt)^r}$, for any $r$, but this agrees with the restriction of $\phi^{\tilt*}\sO_{\bP^{n,\perf}_{K^\tilt}}(1)$ to $V_{(\varpi^\tilt)^r}^\tilt$ for any $r$.  Since these sets cover $X^\tilt$, we are done.

In summary, we have proved the following proposition.
\begin{Proposition}\label{addingasection}
Let $\left(\sL_i,s_j^{(i)},\alpha_i\right)\in\fL_n(X)$, correspond to a map $\phi:X\to\bP^{n,\perf}_{K}$.  Suppose $t_i\in\Gamma(X,\sL_i)$ is a global section such that $\alpha_i\left(t_{i+1}^{\otimes p}\right) = t_i$, and let $\psi:X\to\bP^{n+1}_K$ be the map associated to $\left(\sL_i,\{s_j^{(i)},t_i\},\alpha_i\right)\in\fL_{n+1}(X)$.  Then
\[\phi^{\tilt*}\sO_{\bP^{n,\perf}_{K^\tilt}}(1)\cong\psi^{\tilt*}\sO_{\bP^{n+1,\perf}_{K^\tilt}}(1).\]
\end{Proposition}
Adding sections one at a time by induction completes the proof of Proposition \ref{matchingisos}.

\subsection{Injectivity of $\theta$}\label{injectiveTheta}
With these tools in hand, we can prove the injectivity of $\theta$ for certain perfectoid spaces $X$.  We will first need one more lemma.
\begin{Lemma}\label{GmAction}
Let $\left(\sL_i,s_j^{(i)},\alpha_i\right)\in\fL_n(X)$ correspond to a map $\phi:X\to\bP^{n,\perf}_{K}$.  Fix $\lambda = (\lambda_0,\lambda_1,\cdots)\in\Gamma(X,\sO_X^{\tilt*})$, that is, $\lambda_{i+1}^p = \lambda_i$, so that $\left(\sL_i,\lambda_is_{0}^{(i)},\lambda_i\alpha_i\right)\in\fL_n(X)$ corresponds to a map $\psi:X\to\bP^n_K$.  Then $\phi=\psi$.
\end{Lemma}
\begin{Proof}
Multiplication by $\lambda_i$ for each $i$ gives us an isomorphism $\left(\sL_i,s_j^{(i)},\alpha_i\right)\longtoo{\sim}\left(\sL_i,\lambda_is_{0}^{(i)},\lambda_i\alpha_i\right)$ in $\fL_n(X)$.  Then we are done by Theorem \ref{naturaliso}
\end{Proof}

Before we state the main theorem we make the following definition.
\begin{Definition}\label{walb}
A line bundle $\sL$ on a perfectoid space $X$ is said to be \emph{weakly ample} if for any other line bundle $\sM$, there is some $N>>0$ such that for all $r>N$ we have $\sM\otimes\sL^r$ globally generated.
\end{Definition}
\begin{Theorem}\label{thetaInjects!}
Suppose $X$ is a perfectoid space over $K$. Suppose that $X$ has a weakly ample line bundle and that $H^0(X_{\overline K},\sO_{X_{\overline K}}) = \overline K$, where $\overline K$ is a fixed algebraic closure of $K$.  Then
\[\theta:\Pic X^\tilt\into\ilim_{\sL\mapsto\sL^p}\Pic X.\]
In particular, if $\Pic X$ has no $p$ torsion, then
\[\theta_0:\Pic X^\flat\into\Pic X.\]
\end{Theorem}
\begin{Proof}
Fix $\sL,\sM\in\Pic X^\tilt$ with $\theta(\sL)=\theta(\sM)$.  We first reduce to the case that $\sL,\sM$ are globally generated.  Indeed, letting $\sF$ be a weakly ample line bundle, we have $\theta(\sL\otimes\sF^N) = \theta(\sM\otimes\sF^N)$.  If the result holds for globally generated line bundles, for large enough $N$ we conclude that $\sL\otimes\sF^N\cong\sM\otimes\sF^N$ so that $\sL\cong\sM$.

Next we prove it for the case where $K$ contains all $p$th power roots for all its elements.  Choose generating sections $s_0,\cdots,s_n$ for $\sL$ and $t_0,\cdots,t_r$ of $\sM$, which give us maps
\[\phi^\tilt:X^\tilt\longto\bP^{n,\perf}_{K^\tilt},\]
and
\[\psi^\tilt:X^\tilt\longto\bP^{r,\perf}_{K^\tilt},\]
respectively.  These untilt to
\[\phi:X\longto\bP^{n,\perf}_K,\]
and
\[\psi:X\longto\bP^{r,\perf}_K,\]
which in turn correspond to tuples $\left(\sL_i,s_j^{(i)},\alpha_i\right)\in\fL_n(X)$ and $\left(\sL_i,t_j^{(i)},\beta_i\right)\in\fL_r(X)$.  Notice that $\alpha_i$ and $\beta_i$ differ by an element
\[\lambda_i\in\operatorname{Isom}(\sL_i,\sL_i) = \Gamma(X,\sO_X^*) = K^*.\]
That is, $\alpha_i = \lambda_i\beta_i$.  Choose $p$th power roots $\lambda_i^{1/p^j}$ for each $i,j$ (these exist by assumption), and for all $j$ define:
\begin{eqnarray*}
\tilde t_{j}^{(0)} &=& t_{j}^{(0)}\\
\tilde t_{j}^{(1)} &=& \lambda_0^{-1/p}t_{j}^{(1)}\\
\tilde t_j^{(2)} &=& \lambda_1^{-1/p}\lambda_0^{-1/p^2}t_j^{(2)}\\
{}&\vdots&{}\\
\tilde t_{j}^{(i+1)} &=& \lambda_i^{-1/p}\lambda_{i-1}^{-1/p^2}\cdots\lambda_0^{-1/p^{i+1}}t_{j}^{(i+1)}\\
{}&\vdots&{}
\end{eqnarray*}
Then
\begin{eqnarray*}
\alpha_i\left(\left(\tilde t_j^{(i+1)}\right)^{\otimes p}\right) &=& \lambda_i\beta_i\left(\left(\lambda_{i}^{-1/p}\lambda_{i-1}^{-1/p^2}\cdots\lambda_0^{-1/p^{i+1}}t_{j}^{(i+1)}\right)^{\otimes p}\right)\\
&=&\lambda_i\lambda_i^{-1}\lambda_{i-1}^{1/p}\cdots\lambda_0^{1/p^i}\beta\left(\left(t_j^{(i+1)}\right)^{\otimes p}\right)\\
&=&\lambda_{i-1}^{1/p}\cdots\lambda_0^{1/p^i}t_j^{(i)}\\
&=&\tilde t_j^{(i)}.
\end{eqnarray*}
Therefore the tuple $(\sL_i,\tilde t_j^{(i)},\alpha_i)\in\fL_n(X)$, and it also corresponds to $\psi$ by Lemma \ref{GmAction}.  Furthermore, the isomorphisms corresponding to this data are now $\alpha_i$ in both cases, so that by Proposition \ref{matchingisos}
\[\sL = \phi^{\tilt*}\sO_{\bP^{n,\perf}_{K^\tilt}}(1)\cong\psi^{\tilt*}\sO_{\bP^{n,\perf}_{K^\tilt}}(1) = \sM.\]
For the general case, we let $L/K$ be the extension given by adjoining all $p$th power roots of all elements of $K$.  We have the following diagram.
\begin{cd}
\Pic X_L^\tilt\ar[hookrightarrow,"\theta_L"]{r} & \ilim\Pic X_L\\
\Pic X^\tilt \ar{u}\ar["\theta"]{r} & \ilim\Pic X \uar
\end{cd}
$\theta_L$ injects by the argument we just made.  Furthermore, since $X^\tilt_L\to X^\tilt$ is a pro-\'etale cover of $p$th power degree, the kernel of $\Pic X^\tilt\to\Pic X_L^\tilt$ is $p$th power torsion.  Since $X^\tilt$ is perfect, $\Pic X^\tilt$ has no $p$th power torsion, so the map injects.  Therefore $\theta$ injects.
\end{Proof}
\begin{Open}
In which contexts does $\theta$ surject?
\end{Open}
\bibliography{bib}{}
\bibliographystyle{acm}
\end{document}